\documentclass[11pt]{article}
\usepackage{amsfonts}
\usepackage{eucal}
\usepackage{amssymb,epsfig,latexsym}
\usepackage{graphicx,amsmath,amssymb,latexsym,psfrag, color}
\usepackage{graphics,graphicx}
\newcommand{\dd}{\mathrm d}
\numberwithin{equation}{section}

\newcommand{\EE}{\mathbb{E}}

\newtheorem{Theorem}{Theorem}[section]
\newtheorem{Proposition}[Theorem]{Proposition}

\newtheorem{Definition}[Theorem]{Definition}
\newtheorem{Corollary}[Theorem]{Corollary}
\newtheorem{Lemma}[Theorem]{Lemma}

\title{\huge Almost Sure Central Limit Theorems in Stochastic Geometry}
\author{Giovanni
Luca Torrisi\thanks{Istituto per le Applicazioni del Calcolo, CNR, Roma, Italy.
e-mail: \tt{giovanniluca.torrisi@cnr.it}} \and   Emilio Leonardi\thanks{Dipartimento di Elettronica,
Politecnico di Torino, Italy. e-mail: \tt{emilio.leonardi@polito.it}}  }
\date{}

\allowdisplaybreaks

\begin{document}

\maketitle

\begin{abstract}
We prove an almost sure central limit theorem on the Poisson space, which is perfectly tailored for  stabilizing functionals emerging in stochastic geometry. As a consequence,
we provide almost sure central limit theorems for $(i)$ the total edge length of the $k$-nearest neighbors random graph, $(ii)$
the clique count in random geometric graphs, $(iii)$ the volume of the set approximation via the Poisson-Voronoi tessellation.
\end{abstract}


\baselineskip0.7cm

\section{Introduction}

Let $\{X_n\}_{n\geq 1}$ be a sequence of real-valued random variables and $Z$ a random variable distributed according to the standard normal law.
We say that $\{X_n\}_{n\geq 1}$ satisfies the almost sure central limit theorem (ASCLT) if
\begin{equation}\label{eq:ASCLT}
\lim_{n\to\infty}\frac{1}{\log n}\sum_{k=1}^{n}\frac{1}{k}f(X_k)=\mathbb E[f(Z)],\quad\text{$\forall$ $f\in\mathcal{C}_b(\mathbb R)$, almost surely}
\end{equation}
i.e., the sequence of random measures
\[
\left\{\frac{1}{\log n}\sum_{k=1}^{n}\frac{1}{k}\varepsilon_{X_k}\right\}_{n\geq 2}
\]
converges weakly to the standard normal law, almost surely. Roughly speaking, if an ASCLT holds, then
the Gaussian asymptotic behavior can be observed along individual trajectories of the process.
Here $\mathcal{C}_b(\mathbb R)$ denotes the family of bounded and continuous functions
from $\mathbb R$ to $\mathbb R$ and $\varepsilon_x$ denotes the Dirac measure at $x$.

Property \eqref{eq:ASCLT} should be compared with the classical notion of weak convergence
of $\{X_n\}_{n\geq 1}$ to the random variable $Z$, which can be stated as follows:
\begin{equation}\label{eq:weakconv}
\lim_{n\to\infty}\mathbb{E}[f(X_n)]=\mathbb E[f(Z)],\quad\text{$\forall$ $f\in\mathcal{C}_b(\mathbb R)$.}
\end{equation}
In \eqref{eq:ASCLT} expectation, $\mathbb E$,  is replaced by  logarithmic average, $\frac{1}{\log n}\sum_{k=1}^{n}\frac{1}{k}$, while  almost sure convergence is considered.
Moreover, as it has been recently pointed out in \cite{AN}, a simple application of Skorohod's representation theorem shows that,
in contrast to \eqref{eq:weakconv}, the validity of \eqref{eq:ASCLT} depends  on the whole sequence $\{X_n\}_{n\geq 1}$, and that
\eqref{eq:weakconv} does not imply \eqref{eq:ASCLT}.

The ASCLT for the sequence $\{S_n\}_{n\geq 1}$, where $S_n:=n^{-1/2}\sum_{k=1}^{n}X_k$ and $X_n$, $n\geq 1$, are
independent and identically distributed real-valued random variables with $\EE[X_1]=0$ and $\EE[X_1^2]=1$, was conjectured by L\'evy \cite{L} and proved
in \cite{Br} and \cite{Sc} independently. Since then ASCLTs have attracted a significant interest in the scientific community. For instance, it was proved in \cite{BC}
that (basically) whenever a sequence of independent and identically distributed random variables converges in distribution,
the corresponding almost sure limit theorem holds. The situation is much more complicated for dependent sequences, see e.g. \textcolor{black}{\cite{Ch,Pel,Yo}}.
More recently, in \cite{BNT} the authors  proved ASCLTs for sequences of functionals of general Gaussian fields, with applications to stationary
Gaussian sequences. The main idea was to employ the Malliavin calculus on the Wiener space in order to check the conditions of the so-called Ibragimov-Lifshits criterion \cite{IL}.
A similar approach was followed in \cite{Zh} to provide an ASCLT for sequences of random variables belonging to a fixed Rademacher chaos.

In this paper we prove an ASCLT on the Poisson space, which can be successfully employed to
 stabilizing (or localized) functionals emerging in stochastic geometry. From the point of view of applications,
we provide ASCLTs for the total edge length of the (undirected) $k$-nearest neighbors random graph with $m$th power weighted edges,
the clique count in random geometric graphs, the volume of the set approximation via the Poisson-Voronoi tessellation.

In broad terms, our approach follows the general scheme of \cite{BNT}. We show a (suitable) abstract ASCLT on the Poisson space. To this aim, we exploit some inequalities
for Malliavin's operators on the Poisson space, recently proved in \cite{LPS}, which allow us to express the conditions of the Ibragimov-Lifshits criterion in terms of conditions
involving only the first two Malliavin's derivatives of functionals of the Poisson measure.
The application to stabilizing functionals relies, instead, on the use of some inequalities on moments
of Malliavin's gradients and related probabilities for such functionals, which have recently been proved
\textcolor{black}{in} \cite{LSY}. The specific application to the above mentioned geometric quantities is based on
estimates for the variance of the corresponding functionals proved in \cite{PY,PY2,S,TY} and some
arguments in \cite{LSY} used to check locality of the corresponding functionals.

The paper is organized as follows. In Section \ref{sec:main} we state the main results of the paper, i.e.,
$(i)$ an ASCLT for localized functionals of the Poisson measure, $(ii)$ ASCLTs for the geometric quantities described above,
$(iii)$ an abstract ASCLT on the Poisson space.
The proofs of such results are given in Sections \ref{sec:stabilizing}, \ref{sec:rg} and \ref{sec:poissonspace}, respectively.

\section{Main results}\label{sec:main}

Throughout this paper we consider a probability space $(\Omega,\mathcal F,\mathbb P)$,
let $(\mathbb X,\mathcal X)$ be an arbitrary measurable space, let $\mu$ be a $\sigma$-finite measure on $\mathbb X$ with $\mu(\mathbb X)>0$
and denote by $\eta$ the Poisson random measure on $\mathbb X$ with intensity measure $\mu$.
Formally, we view $\eta$ as a random element on $\bold N_{\mathbb X}$, the space of integer-valued and simple $\sigma$-finite measures
$\nu$ equipped with the smallest $\sigma$-field $\mathcal{N}_{\mathbb X}$ that makes the mappings $\nu\mapsto\nu(B)$ measurable for any $B\in\mathcal X$.
When explicitly stated, for ease of notation, we identify a counting measure with its support.

\subsection{An ASCLT for stabilizing functionals on the Poisson space}

Let $\mathbb Y$ be a full-dimensional Borel set of $\mathbb R^d$ and let
$\eta$ be the Poisson random measure on $\mathbb X:=\mathbb{R}^d$
with intensity measure $\mu(\mathrm{d}x):=\ell_d(\mathrm{d}x)$, the Lebesgue measure on $\mathbb R^d$.
Here and in the next subsection, we identify a counting measure $\nu\in\bold{N}_{\mathbb R^d}$ with its support and we consider the statistic
\begin{equation}\label{eq:stabfunc}
H_n:=\sum_{x\in\eta_n\cap\mathbb Y}\xi_n(x,\eta_n\cap\mathbb Y),\quad n\geq 1
\end{equation}
where $\eta_n:=n^{-1/d}\eta$, i.e., $\eta_n$ is the Poisson random measure on $\mathbb R^d$ with intensity measure $n\ell_d(\dd x)$.
The so-called score functions $\{\xi_n\}_{n\geq 1}$ are measurable functions from $\mathbb Y\times\bold N_\mathbb Y$ to $\mathbb R$ and represent the local contribution to the global statistic $H_n$.
In order to introduce the notion of stabilizing functional considered in this paper,
we need some definitions.

For $n\geq 1$, a measurable map $R_n:\mathbb Y\times\bold N_\mathbb Y\to [0,\infty)$ is called a $\it{radius\,\,of\,\,stabilization}$ for the score function $\xi_n$ if,
for all $(x,\nu)\in \mathbb Y\times\bold N_\mathbb Y$
and any finite subset $\mathcal A\subset \mathbb Y$ with at most $7$ points,
\begin{equation}\label{eq:RadStab}
\xi_n(x,(\nu\cup\{x\}\cup\mathcal A)\cap B(x,R_n(x,\nu\cup\{x\})))=\xi_n(x,\nu\cup\{x\}\cup\mathcal A).
\end{equation}
Here $B(x,r):=\{y\in\mathbb Y:\,\,\|x-y\|_d\leq r\}$ denotes the closed ball of radius $r\geq 0$ centered at $x\in\mathbb Y$
and $\|\cdot\|_d$ is the Euclidean norm on $\mathbb R^d$. Loosely speaking, the notion of radius of stabilization
says that the value of $\xi_n$ at $x$ is wholly determined by the points of $\nu$ at distance at most $R_n(x,\nu\cup\{x\})$ from $x$.
\textcolor{black}{We emphasize that the assumption \lq\lq $\mathcal{A}$ has at most $7$ points\rq\rq\, is not required according to  the classical definition of
radius of stabilization. In this paper we need this extra hypothesis to place ourselves in the same framework of \cite{LSY}.
A similar comment applies to the following definitions of score functions satisfying a moment condition and score functions with
exponentially fast decay.}

The score functions $\{\xi_n\}_{n\geq 1}$ are called $\it{exponentially\,\,stabilizing}$ if there exist radii of stabilization
$\{R_n\}_{n\geq 1}$ and constants $C_{stab}$, $c_{stab}$, $\alpha_{stab}\in\textcolor{black}{\mathbb R_+:=(0,\infty)}$ such that, for $x\in \mathbb Y$, $r\geq 0$ and $n\geq 1$,
\begin{equation}\label{eq:expstab}
\mathbb{P}(R_n(x,(\eta_n\cap \mathbb Y)\cup\{x\})\geq r)\leq C_{stab}\exp(-c_{stab}n^{\alpha_{stab}/d}r^{\alpha_{stab}}).
\end{equation}

Let $p'\in [0,\infty)$ be given. We say that the score functions $\{\xi_n\}_{n\geq 1}$ satisfy a $(4+p')$$\it{th}$ $\it{moment\,\,condition}$ if there
is a constant $C_{p'}\in\mathbb{R}_+$ such that for any finite subset $\mathcal{A}\subset\mathbb{Y}$ with at most $7$ points,
\begin{equation}\label{eq:boundedmoment}
\sup_{n\geq 1}\sup_{x\in\mathbb Y}\mathbb{E}[|\xi_n(x,(\eta_n\cap \mathbb Y)\cup\{x\}\cup\mathcal A)|^{4+p'}]\leq C_{p'}.
\end{equation}

Let $K$ be a Borel subset of $\mathbb Y$ and put $\mathrm{d}(x,K):=\inf_{y\in K}\|x-y\|_d$, $x\in\mathbb Y$. We say that the score functions $\{\xi_n\}_{n\geq 1}$
$\it{decay\,\,exponentially\,\,fast\,\,with\,\,the\,\,distance\,\,to\,\,K}$ if there are constants $C_K,c_K,\alpha_K\in\mathbb{R}_+$
such that for any finite subset $\mathcal{A}\subset\mathbb Y$ with at most $7$ points,
\begin{equation}\label{eq:expdecay}
\mathbb{P}(\xi_n(x,(\eta_n\cap\mathbb Y)\cup\{x\}\cup\mathcal A)\neq 0)\leq C_K\exp(-c_K n^{\alpha_K/d}\mathrm{d}(x,K)^{\alpha_K}),
\end{equation}
for $x\in\mathbb Y$ and $n\geq 1$.

\begin{Definition}
We say that the functionals $H_n$, $n\geq 1$, defined by \eqref{eq:stabfunc} are stabilizing if
their score functions enjoy the properties \eqref{eq:expstab}, \eqref{eq:boundedmoment} for some $p'\in (0,1]$, and \eqref{eq:expdecay}
\textcolor{black}{for some Borel set $K\subseteq\mathbb Y$.}
\end{Definition}

Hereon, we consider stabilizing functionals $H_n$ and, for a fixed ${p''}\in (0,p')$, the quantity
\begin{equation}\label{eq:IK}
I_{K,n}:=n\int_{\mathbb Y}\exp\left(-\frac{{p''}\,c_{p'}n^{\alpha/d}\mathrm{d}(x,K)^\alpha}{2^{2\alpha+3}(4+{p''})}\right)\,\dd x,\quad n\geq 1
\end{equation}
where
\begin{equation}\label{eq:cpalpha}
c_{p'}:=p'(c_{stab}\wedge c_K),\quad\alpha:=\alpha_{stab}\wedge\alpha_K
\end{equation}
and we denote by $a\wedge b$ the minimum between $a,b\in\mathbb R$. \textcolor{black}{Throughout this paper, we use the standard Landau notation, i.e.,
given two sequences $\{a_n\}_{n\geq 1},\{b_n\}_{n\geq 1}\subset\mathbb R_+$, we write $a_n=O(b_n)$ if $\limsup_{n\to\infty}\frac{a_n}{b_n}<\infty$.}

The following theorem holds.

\begin{Theorem}\label{thm:ASCLTstabilizing}
If the functionals $H_n$, $n\geq 1$, defined by \eqref{eq:stabfunc} are stabilizing, $n^\tau=O(\mathbb{V}ar(H_n))$ and $I_{K,n}=O(n^\tau)$, for some $\tau\in (0,2)$, then
$\{F_n\}_{n\geq 1}$, $F_n:=(H_n-\mathbb{E}[H_n])/\sqrt{\mathbb{V}ar(H_n)}$, satisfies the ASCLT.
\end{Theorem}

\textcolor{black}{We remark that in the statement of Theorem \ref{thm:ASCLTstabilizing} the quantities $\mathbb{V}ar(H_n)$ and $I_{K,n}$ are finite and strictly positive
for any $n\geq 1$ (this is an implicit consequence of the Landau notation, as defined above, which is employed to compare sequences in $\mathbb R_+$). We also emphasize that,
although Theorem \ref{thm:ASCLTstabilizing} holds for any scaling regime $\tau\in(0,2)$, all the examples considered in this paper fall either in the class of
the volume order scaling (i.e., $\tau=1$) or in the class of the surface area order scaling (i.e. $\tau=1-1/d$).}

\subsection{ASCLTs for stochastic geometry models}

Theorem \ref{thm:ASCLTstabilizing} finds a natural application to models emerging in stochastic geometry, such as the total edge length of the $k$-nearest neighbors random graph,
the clique count in random geometric graphs and the volume of the set approximation via the Poisson-Voronoi tessellation.
%The statements of the corresponding results are given in the next three , respectively.

\subsubsection{An ASCLT for the total edge length of the $k$-nearest neighbors random graph}

Let $\mathbb{Y}$ be a full-dimensional Borel set of $\mathbb R^d$ and $n,k\geq 1$ two integers. For $x\in\eta_n\cap\mathbb Y$, we denote by $V_k(x,\eta_n\cap\mathbb Y)$
the set of the $k$ nearest neighbors of $x$, i.e., the $k$ closest points of $x$ in $(\eta_n\cap\mathbb Y)\setminus\{x\}$. The (undirected) $k$-nearest neighbors random graph
$NG_k(\eta_n\cap\mathbb Y)$ is the random graph with vertex set $\eta_n\cap\mathbb Y$ obtained by including an edge $\{x,y\}$ if $y\in V_k(x,\eta_n\cap \mathbb Y)$ and/or
$x\in V_k(y,\eta_n\cap \mathbb Y)$. For all $m\geq 0$, we set
\[
\xi^{(m)}(x,\eta_n\cap\mathbb Y):=\sum_{y\in V_k(x,\eta_n\cap\mathbb Y)}\rho^{(m)}(x,y),
\]
where
\begin{align}
\rho^{(m)}(x,y):&=\bold{1}\{x\in V_k(y,\eta_n\cap\mathbb Y)\}\frac{\|x-y\|_d^m}{2}
%\nonumber\\
%&\qquad\qquad
%+\bold{1}\{x\in V_k(y,\eta_n\cap\mathbb Y),y\notin V_k(x,\eta_n\cap\mathbb Y)\}\|x-y\|_d^q\nonumber\\
%&\qquad\qquad
+\bold{1}\{x\notin V_k(y,\eta_n\cap\mathbb Y)\}\|x-y\|_d^m.\nonumber
\end{align}
The total edge length of the (undirected) $k$-nearest neighbors random graph on $\eta_n\cap\mathbb Y$ with $m$th power weighted edges is
\[
L_{NG_k}^{(m)}(\eta_n\cap\mathbb Y):=\sum_{x\in\eta_n\cap\mathbb Y}\xi^{(m)}(x,\eta_n\cap\mathbb Y).
\]
Note that $L_{NG_k}^{(0)}(\eta_n\cap\mathbb Y)$ is the number of edges and $L_{NG_k}^{(1)}(\eta_n\cap\mathbb Y)$ is the total edge length. When $\mathbb Y$
is a full-dimensional compact and convex subset of $\mathbb R^d$, the central limit theorem
and a quantitative central limit theorem (in the Kolmogorov distance) for the sequence $\{L_{NG_k}^{(m)}(\eta_n\cap\mathbb Y)\}_{n\geq 1}$ are proved
in \cite{PY} and \cite{LPS}, respectively. Here, we provide an ASCLT for such sequence.
%$\{L_{NG_k}^{(q)}(\eta_n\cap\mathbb Y)\}_{n\geq 1}$.

Set
\begin{equation}\label{eq:totaledge}
F_n:=
%\frac{H_n-\mathbb{E}[H_n]}{\sqrt{\mathbb{V}ar(H_n)}}=
\frac{L_{NG_k}^{(m)}(\eta_n\cap\mathbb Y)-\mathbb{E}[L_{NG_k}^{(m)}(\eta_n\cap\mathbb Y)]}
{\sqrt{\mathbb{V}ar(L_{NG_k}^{(m)}(\eta_n\cap\mathbb Y))}},\quad n\geq 1.
\end{equation}
The following corollary holds.

\begin{Corollary}\label{cor:totaledge}
If $\mathbb Y$ is a full-dimensional compact and convex subset of $\mathbb R^d$, then $\{F_n\}_{n\geq 1}$ defined by \eqref{eq:totaledge} satisfies the ASCLT.
\end{Corollary}

\subsubsection{An ASCLT for the clique count in random geometric graphs}

Let $\mathbb{Y}$ be a full-dimensional subset of $\mathbb{R}^d$ and consider the random geometric graph $G(\eta_n\cap\mathbb Y,r)$, $n\geq 1$, $r\in\mathbb{R}_+$, where two nodes $x,y\in\eta_n\cap\mathbb Y$
are joined with an edge if $\|x-y\|_d\leq r$. We recall that $k+1$ nodes of $G(\eta_n\cap\mathbb Y,r)$ form a clique of order $k+1$ if each pair of them is connected by an edge. The number of cliques
of order $k+1$ in $G(\eta_n\cap\mathbb Y,r)$, denoted by $\mathcal{C}_k(\eta_n\cap\mathbb Y,r)$, is a central statistic in topological data analysis.
%When $\mathbb Y$ is a
%full-dimensional subset of $\mathbb R^d$
A quantitative central limit theorem (in the Kolmogorov distance) for the sequence $\{\mathcal{C}_k(\eta_n\cap\mathbb Y,rn^{-1/d})\}_{n\geq 1}$ has been proved in \cite{Penrosebook},
\textcolor{black}{Chapter 3.} Here, we give an ASCLT for such sequence.

Set
\begin{equation}\label{eq:cliquecount}
F_n:=\frac{\mathcal{C}_k(\eta_n\cap\mathbb Y,rn^{-1/d})-\mathbb{E}[\mathcal{C}_k(\eta_n\cap\mathbb Y,rn^{-1/d})]}{\sqrt{\mathbb{V}ar(\mathcal{C}_k(\eta_n\cap\mathbb Y,rn^{-1/d}))}},\quad
n\geq 1.
\end{equation}
The following corollary holds.

\begin{Corollary}\label{cor:cliquecount}
If $\mathbb Y$ is a full-dimensional subset of $\mathbb R^d$ such that $\ell_d(\mathbb Y)<\infty$, then $\{F_n\}_{n\geq 1}$ defined by \eqref{eq:cliquecount} satisfies the ASCLT.
\end{Corollary}

\subsubsection{An ASCLT for the volume of the set approximation via the Poisson-Voronoi tessellation}

Let $\mathbb Y:=[-1/2,1/2]^d$, $d\geq 2$, and let $A\subset (-1/2,1/2)^d$ be a full-dimensional subset of $\mathbb R^d$.
%be a full-dimensional subset of $\mathbb R^d$.
Assume that we observe $\eta_n\cap\mathbb Y$ and that the only information about $A$ at our disposal is which points of
$\eta_n\cap\mathbb Y$ lie in $A$, i.e., we have the partition of the process $\eta_n\cap\mathbb Y$ into $\eta_n\cap A$ and $(\eta_n\cap\mathbb Y)\setminus A$.
To reconstruct the set $A$ just by the information at our disposal, we approximate $A$ by its Poisson-Voronoi approximation $A_n$, i.e.,
the set of all points in $\mathbb Y$ which are closer to $\eta_n\cap A$ than to $(\eta_n\cap\mathbb Y)\setminus A$. Formally,
\[
A_n:=\bigcup_{x\in\eta_n\cap A}\mathrm{C}(x,\eta_n\cap\mathbb Y),
\]
where $\mathrm{C}(x,\eta_n\cap\mathbb Y)$ is the Poisson-Voronoi cell generated by $\eta_n\cap\mathbb Y$ with nucleus $x\in\eta_n\cap\mathbb Y$, i.e.,
\[
\mathrm{C}(x,\eta_n\cap\mathbb Y):=\{y\in\mathbb Y:\,\,\|y-x\|_d\leq\|y-z\|_d\,\,\forall\,\,z\in\eta_n\cap\mathbb Y\}.
\]
Note that $\mathrm{C}(x,\eta_n\cap\mathbb Y)$ is a random convex polytope and that the Poisson-Voronoi tessellation (or mosaic) of $\mathbb Y$, i.e.,
the family $\{\mathrm{C}(x,\eta_n\cap\mathbb Y)\}_{x\in\eta_n\cap\mathbb Y}$, is a partition of $\mathbb Y$.

If $A$ is compact and convex it follows, respectively from \cite{P} and
\cite{R}, that, as $n\to\infty$, $\ell_d(A_n)\to\ell_d(A)$ almost surely, and $\mathbb{E}[\ell_d(A_n)]\to\ell_d(A)$.
%\quad $\ell_d(A\Delta A_n)\to 0$,
%\quad as $n\to\infty$, $\mathbb{P}$-a.s.}
%\end{equation*}
%and $\mathcal{H}_{d-1}(\partial A_n)$ approximate in $L^2$ a scalar multiple of $\mathcal{H}_{d-1}(\partial A)$, as $n\to\infty$.
%Here, $A_n:=A(\eta_n\cap\mathbb Y)$, $\Delta$ is the operation of symmetric difference of sets and
%$\mathcal{H}_{d-1}$ denotes the $(d-1)$-dimensional Hausdorff measure.
A quantitative central limit theorem (in the Kolmogorov distance) for the sequence $\{\ell_d(A_n)\}_{n\geq 1}$
%, $\ell_d(A\Delta A_n)$ and $\mathcal{H}_{d-1}(\partial A_n)$
is proved in \cite{LSY}. Here, we give an ASCLT for such sequence.

Put
\begin{equation}\label{eq:voronoi}
F_n:=\frac{\ell_d(A_n)-\mathbb E[\ell_d(A_n)]}{\sqrt{\mathbb{V}ar(\ell_d(A_n))}},\quad n\geq 1.
\end{equation}
%where either $\varphi(\cdot):=\ell_d(\cdot)$ or $\varphi(\cdot):=\ell_d(A\Delta\cdot)$ or $\varphi(\cdot):=\mathcal{H}_{d-1}(\partial\,\cdot)$.
The following corollary holds.

\begin{Corollary}\label{cor:voronoimosaic}
If $A$ is compact and convex
%or $A$ is closed and $\partial A$ contains a $(d-1)$-dimensional $C^2$-submanifold,
then $\{F_n\}_{n\geq 1}$ defined by \eqref{eq:voronoi} satisfies the ASCLT.
\end{Corollary}

\subsection{An ASCLT on the Poisson space}

As already mentioned in the Introduction, in this section we present an abstract
ASCLT for the Poisson space, which provides theoretical foundation for the  proof of Theorem \ref{thm:ASCLTstabilizing}. To state such general result we need to introduce some additional notation.

By $L_{\eta}^r$, $r\in\mathbb R_+$, we denote the space of all random variables $F\in L^r(\mathbb P)$ such that $F=f(\eta)$ $\mathbb P$-a.s. for some measurable function $f:\bold{N}_{\mathbb X}\to\mathbb R$,
where $L^r(\mathbb P)$ is the set of random variables $X:\Omega\to\mathbb R$ such that $\EE[|X|^r]<\infty$.
For $F\in L_{\eta}^2$, $F=f(\eta)$, and $x_1,x_2\in\mathbb X$ we define the Malliavin derivative operators
\[
D_{x_1}F:=f(\eta+\varepsilon_{x_1})-f(\eta),
\]
\[
D_{x_1,x_2}^{2}F:=D_{x_2}(D_{x_1}F)=f(\eta+\varepsilon_{x_1}+\varepsilon_{x_2})-f(\eta+\varepsilon_{x_2})-f(\eta+\varepsilon_{x_1})+f(\eta).
\]

\textcolor{black}{For $H_i\in L_\eta^2$ with $\mathbb{V}ar(H_i)>0$, $i=1,2,3,4$,} and positive constants $p,q\in\mathbb R_+$ (defined in the statement of Theorem \ref{cor:ASCLTAPP}), we put
\begin{align}
&\bold{\Gamma}(H_1,H_2,H_3,H_4)^2\nonumber\\
&\qquad
:=\bold{\Lambda}(H_1,H_2,H_3,H_4)^2
+\frac{1
%c^{4/(4+q)}
}{\mathbb{V}ar(H_1)\mathbb{V}ar(H_2)}
\int_{\mathbb{X}}\psi_x(H_1,q/(2(4+q)))\psi_x(H_2,q/(2(4+q)))\mu(\dd x),\nonumber
\end{align}
\[
\psi_x(G,\beta):=\int_{\mathbb{X}}\mathbb{P}(D_{x_1,x}^{2}G\neq 0)^{\beta}\mu(\dd x_1),\quad x\in\mathbb X,\,G\in L_\eta^2,\,\beta\in\mathbb{R}_+,
\]
\begin{align}
%&
\bold{\Lambda}(H_1,H_2,H_3,H_4)^2
%\nonumber\\
%&\qquad
:=\frac{1
%4c^{2/(4+q)+2/(4+p)}
}
{\prod_{i=1}^{4}\mathbb{V}ar(H_i)^{1/2}}
\int_{\mathbb{X}}\psi_x(H_1,q/(4(4+q)))\psi_x(H_2,q/(4(4+q)))
\mu(\dd x),\nonumber
\end{align}
\[
\bold{\Theta}(H_1,H_2):=\frac{1
%c^{2/(4+p)}
}{\prod_{i=1}^{2}\mathbb{V}ar(H_i)^{1/2}}
\int_{\mathbb X}\mathbb{P}(D_x H_1\neq 0)^{p/(4(4+p))}\mathbb{P}(D_x H_2\neq 0)^{p/(4(4+p))}\mu(\mathrm{d}x),
\]
\begin{align*}
\bold{\Gamma}_1(H_1)^2:=\bold{\Gamma}(H_1,H_1,H_1,H_1)^2,
%\frac{4c^{4/(4+q)}}{\mathbb{V}ar(H_1)^2}\int_{\mathbb{X}}\psi_x(H_1,q/(4(4+q)))^2\mu(\dd x)+\frac{c^{4/(4+q)}}{\mathbb{V}ar(H_1)^2}\int_{\mathbb{X}}
%\psi_x(H_1,q/(2(4+q)))^2\mu(\dd x),\nonumber
\end{align*}
\[
\bold{\Gamma}_2(H_1):=\frac{1
%2^{-1/2}c^{3/(4+p)}
}{\mathbb{V}ar(H_1)^{3/2}}\int_{\mathbb{X}}\mathbb{P}(D_x H_1\neq 0)^{(1+p)/(4+p)}\mu(\dd x).
\]
%Note that $\bold{\Gamma}_1(H)^2=\bold{\Gamma}(H,H,H,H)^2$ and
%\begin{equation*}
%\bold{\Lambda}(H_1,H_2,H_3,H_4)=\bold{\Lambda}(H_2,H_1,H_4,H_3)\quad\text{and}\quad\bold{\Gamma}(H_1,H_2,H_3,H_4)=\bold{\Gamma}(H_2,H_1,H_4,H_3).
%\end{equation*}
The following theorem holds, where we denote by $a\vee b$ the maximum between $a,b\in\mathbb R$.

\begin{Theorem}\label{cor:ASCLTAPP}
Assume $\{H_n\}_{n\geq 1}\subset L_\eta^2$, $\mathbb{V}ar(H_n)>0$, $n\geq 1$, and that there exist constants $c,p,q\in\mathbb{R}_+$ such that
\begin{equation}\label{eq:os1app}
\sup_{n\geq 1}\mathbb{E}[|D_{x_1}H_n|^{4+p}]\leq c,
\quad\text{$\mu$-a.e. $x_1\in\mathbb X$}
\end{equation}
\begin{equation}\label{eq:os2app}
\sup_{n\geq 1}\mathbb{E}[|D_{x_1,x_2}^{2}H_n|^{4+q}]\}\leq c,
\quad\text{$\mu^{\otimes 2}$-a.e. $(x_1,x_2)\in\mathbb X^{2}$.}
\end{equation}
Moreover, suppose
\begin{equation}\label{hyp:convgammaH}
\lim_{n\to\infty}\bold{\Gamma}_1(H_n)=\lim_{n\to\infty}\bold{\Gamma}_2(H_n)=0,
\end{equation}
\begin{equation}\label{hyp:series1convH}
\sum_{n\geq 2}\frac{1}{n(\log n)^2}\sum_{k=1}^{n}\frac{\bold{\Gamma}_i(H_k)}{k}<\infty,\quad\text{$i=1,2$}
\end{equation}
\begin{equation}\label{hyp:series2convH}
\sum_{n\geq 2}\frac{1}{n(\log n)^3}\sum_{l=1}^{n}\sum_{k=1}^{l}\frac{\bold{\Theta}(H_k,H_l)}{kl}<\infty,
\end{equation}
\begin{equation}\label{hyp:series3convH}
\sum_{n\geq 2}\frac{1}{n(\log n)^3}\sum_{l=1}^{n}\sum_{k=1}^{l}\frac{\bold{\Gamma}(H_k,H_l,H_k,H_l)}{kl}<\infty,
\end{equation}
\begin{equation}\label{hyp:series4convH}
\sum_{n\geq 2}\frac{1}{n(\log n)^3}\sum_{l=1}^{n}\sum_{k=1}^{l}\frac{\bold{\Lambda}_{max}(H_k,H_l)}{kl}<\infty,
\end{equation}
where
\begin{align}
\bold{\Lambda}_{max}(H_k,H_l):=&\bold{\Lambda}(H_k,H_k,H_k,H_l)\vee\bold{\Lambda}(H_l,H_l,H_l,H_k)\vee
\bold{\Lambda}(H_k,H_k,H_l,H_l)\vee\bold{\Lambda}(H_l,H_l,H_k,H_k)\nonumber\\
&\vee\bold{\Lambda}(H_l,H_k,H_k,H_k)\vee\bold{\Lambda}(H_k,H_l,H_l,H_l)
\vee\bold{\Lambda}(H_k,H_l,H_l,H_k).\nonumber
\end{align}
Then $\{F_n\}_{n\geq 1}$, $F_n:=(H_n-\mathbb{E}[H_n])/\sqrt{\mathbb{V}ar(H_n)}$, satisfies the ASCLT.
\end{Theorem}
%This theorem can be successfully applied to  stabilizing functionals that emerge in stochastic geometry.
%The ASCLT corresponding to such functionals is stated in the next subsection.

\section{Proof of Theorem \ref{thm:ASCLTstabilizing}}\label{sec:stabilizing}

The proof of Theorem \ref{thm:ASCLTstabilizing} exploits some recent inequalities, stated in Subsection \ref{sec:prellemmaSTAB},
concerning moments and probabilities of Malliavin's gradients of
localized functionals due to Lachi\'eze-Rey, Schulte and Yukich \cite{LSY}. The proof of the theorem is then given in Subsection \ref{sec:proveThmSTAB}.

We start introducing some notation. Throughout we assume that \textcolor{black}{the functionals $H_n$, $n\geq 1$, defined by \eqref{eq:stabfunc} are stabilizing.}
Recalling that here we identify a counting measure $\nu\in\bold{N}_{\mathbb{R}^d}$ with its support, for $y_1,y_2\in\mathbb Y$, we set
\[
\mathcal{D}_{y_1}H_n:=\sum_{x\in(\eta_n\cap\mathbb Y)\cup\{y_1\}}\xi_n(x,(\eta_n\cap\mathbb Y)\cup\{y_1\})-\sum_{x\in\eta_n\cap\mathbb Y}\xi_n(x,\eta_n\cap\mathbb Y)
\]
and
\[
\mathcal{D}_{y_1,y_2}^{2}H_n:=\mathcal{D}_{y_2}(\mathcal{D}_{y_1}H_n)=\mathcal{D}_{y_1,y_2}^{2,+}H_n-\mathcal{D}_{y_1}H_n,
\]
where
\[
\mathcal{D}_{y_1,y_2}^{2,+}H_n:=\sum_{x\in(\eta_n\cap\mathbb Y)\cup\{y_1,y_2\}}\xi_n(x,(\eta_n\cap\mathbb Y)\cup\{y_1,y_2\})-\sum_{x\in(\eta_n\cap\mathbb Y)\cup\{y_2\}}
\xi_n(x,(\eta_n\cap\mathbb Y)\cup\{y_2\}).\nonumber
\]
We remark that these quantities are well-defined thanks to Lemma \ref{le:LSY1} of the next subsection.

\subsection{Preliminary lemmas}\label{sec:prellemmaSTAB}

The following lemma is an immediate consequence of Lemma 5.5 in \cite{LSY}.

\begin{Lemma}\label{le:LSY1}
\textcolor{black}{If the functionals $H_n$, $n\geq 1$, defined by \eqref{eq:stabfunc} are stabilizing}, then, for any fixed ${p''}\in (0,p')$,
there exists a constant $c({p''})\in\mathbb{R}_+$ $($only depending on the constants $C_{stab}$, $c_{stab}$, $\alpha_{stab}$,
$p'$ and $C_{p'}$ appearing in the definition of stabilizing functional$)$ such that
\[
\sup_{n\geq 1}\sup_{y\in\mathbb Y}\mathbb{E}[|\mathcal{D}_y H_n|^{4+{p''}}]\leq c({p''})\quad\text{and}
\quad
\sup_{n\geq 1}\sup_{y_1,y_2\in\mathbb{Y}}
\mathbb{E}[|\mathcal{D}_{y_1,y_2}^{2,+}H_n|^{4+{p''}}]\leq c({p''}).
\]
\end{Lemma}

The following lemma is an immediate consequence of Lemma 5.10 in \cite{LSY} (see the inequalities $(5.6)$ and $(5.8)$ therein).

\begin{Lemma}\label{le:LSY2}
\textcolor{black}{If the functionals $H_n$, $n\geq 1$, defined by \eqref{eq:stabfunc} are stabilizing}, then, for any fixed
$\beta\in\mathbb{R}_+$ there exists a constant $\tilde{C}_\beta\in\mathbb R_+$ $($not depending on $n$$)$ such that
\[
n\int_{\mathbb{Y}}\Biggl(n\int_{\mathbb Y}\mathbb{P}(
\mathcal{D}_{y_1,y_2}^2H_n \neq 0)^\beta\,\dd y_1\Biggr)^2\,\dd y_2
\leq\tilde{C}_\beta n\int_{\mathbb Y}\exp\left(-\frac{c_{p'}\beta}{2^{2\alpha+1}}(n^{1/d}\dd (x,K))^\alpha\right)\,\dd x
\]
and
\[
n\int_{\mathbb Y}\mathbb{P}(\mathcal{D}_y H_n\neq 0)^\beta\,\dd y\nonumber
\leq\tilde{C}_\beta n\int_{\mathbb Y}\exp\left(-\frac{c_{p'}\beta}{2^{\alpha+1}}(n^{1/d}\dd (x,K))^\alpha\right)\,\dd x,
\]
where the constants $c_{p'}$ and $\alpha$ are defined in \eqref{eq:cpalpha}.
\end{Lemma}

We conclude this subsection with a lemma which relates the operators $D$ and $D^2$ with the operators $\mathcal{D}$
and $\mathcal{D}^2$, respectively.

\begin{Lemma}\label{le:malliavinderivatives}
\textcolor{black}{If the functionals $H_n$, $n\geq 1$, defined by \eqref{eq:stabfunc} are stabilizing} and $I_{K,n}<\infty$, then $H_n\in L^2(\mathbb P)$ and,
for any $x_1,x_2\in\mathbb{R}^d$,
\begin{align}
D_{x_1}H_n=\bold{1}\{n^{-1/d}x_1\in\mathbb Y\}\mathcal{D}_{n^{-1/d}x_1}H_n\label{eq:der1}
\end{align}
and
\begin{align}
D_{x_1,x_2}^2 H_n&=\bold{1}\{n^{-1/d}x_1,n^{-1/d}x_2\in\mathbb Y\}\mathcal{D}_{n^{-1/d}x_1,n^{-1/d}x_2}^2H_n.\label{eq:der2}
\end{align}
\end{Lemma}
\noindent{\it Proof.} We divide the proof in two steps. In the first step we check the square integrability of $H_n$, in the second step
we verify the relations among the operators.\\
\noindent{\it Step\,\,1:\,\,Checking\,\,$H_n\in L^2(\mathbb P)$,\,\,$n\geq 1$.} Throughout this proof ${p''}\in (0,p')$ is the constant
appearing  in the definition of $I_{K,n}$. We shall check later on that, for any $0<\gamma\leq 2^{2\alpha+3}(4+{p''})$ and $x\in\mathbb Y$,
\begin{equation}\label{eq:ineqexp}
\exp(-c_K n^{\alpha_K/d}\mathrm{d}(x,K)^{\alpha_K}/\gamma)\leq\hat{C}\exp(-{p''}\,c_{p'}n^{\alpha/d}\mathrm{d}(x,K)^{\alpha}/(2^{2\alpha+3}(4+{p''}))),
\end{equation}
for a suitable constant $\hat C\in\mathbb R_+$ not depending on $n$,
where the constants $c_{p'}$ and $\alpha$ are defined in \eqref{eq:cpalpha}. Setting $M:=\sum_{x\in\eta_n\cap\mathbb Y}\bold{1}(\xi_n(x,\eta_n\cap\mathbb Y)\neq 0)$, by Jensen's inequality we have
\begin{align}
H_n^2&=M^2\left(\sum_{x\in\eta_n\cap\mathbb Y}\frac{\bold{1}\{\xi_n(x,\eta_n\cap\mathbb Y)\neq 0\}}{M}\xi_n(x,\eta_n\cap\mathbb Y)\right)^2\nonumber\\
&\leq M\sum_{x\in\eta_n\cap\mathbb Y}\bold{1}\{\xi_n(x,\eta_n\cap\mathbb Y)\neq 0\}\xi_n(x,\eta_n\cap\mathbb Y)^2\nonumber\\
&=\sum_{x,x'\in\eta_n\cap\mathbb Y}\bold{1}\{\xi_n(x',\eta_n\cap\mathbb Y)\neq 0,\,\xi_n(x,\eta_n\cap\mathbb Y)\neq 0\}\xi_n(x,\eta_n\cap\mathbb Y)^2\nonumber\\
&=\sum_{x\in\eta_n\cap\mathbb Y}\bold{1}\{\xi_n(x,\eta_n\cap\mathbb Y)\neq 0\}\xi_n(x,\eta_n\cap\mathbb Y)^2\nonumber\\
&\qquad\qquad\qquad
+\sum_{x,x'\in\eta_n\cap\mathbb Y:\,x\neq x'}\bold{1}\{\xi_n(x',\eta_n\cap\mathbb Y)\neq 0,\,\xi_n(x,\eta_n\cap\mathbb Y)\neq 0\}\xi_n(x,\eta_n\cap\mathbb Y)^2.\nonumber
\end{align}
By the (multivariate) Mecke formula (see e.g. \cite{LastPenbook}, \textcolor{black}{Chapter 4}), the Cauchy-Schwarz inequality, H\"older's inequality, \eqref{eq:boundedmoment},
\eqref{eq:expdecay} and \eqref{eq:ineqexp}, we get
\begin{align}
&\mathbb{E}[H_n^2]\leq n\int_{\mathbb Y}\mathbb{E}[\bold{1}\{\xi_n(x,(\eta_n\cap\mathbb Y)\cup\{x\})\neq 0\}\xi_n(x,(\eta_n\cap\mathbb Y)\cup\{x\})^2]\,\dd x\nonumber\\
&\qquad
+n^2\int_{\mathbb{Y}^2}\mathbb{E}[\bold{1}\{\xi_n(x',(\eta_n\cap\mathbb Y)\cup\{x,x'\})\neq 0,\,\xi_n(x,(\eta_n\cap\mathbb Y)\cup\{x,x'\})\neq 0\}\nonumber\\
&\qquad\qquad\qquad\qquad\qquad\qquad\qquad\qquad\qquad
\times\xi_n(x,(\eta_n\cap\mathbb Y)\cup\{x,x'\})^2]\,\dd x\dd x'\nonumber\\
&\leq n\int_{\mathbb Y}\mathbb{P}(\xi_n(x,(\eta_n\cap\mathbb Y)\cup\{x\})\neq 0)^{1/2}\mathbb{E}[\xi_n(x,(\eta_n\cap\mathbb Y)\cup\{x\})^4]^{1/2}\,\dd x\nonumber\\
&\qquad
+n^2\int_{\mathbb{Y}^2}\mathbb{E}[\bold{1}\{\xi_n(x',(\eta_n\cap\mathbb Y)\cup\{x,x'\})\neq 0,\,\xi_n(x,(\eta_n\cap\mathbb Y)\cup\{x,x'\})\neq 0\}]^{1/2}\nonumber\\
&\qquad\qquad\qquad\qquad\qquad\qquad\qquad\qquad\qquad
\times\mathbb{E}[\xi_n(x,(\eta_n\cap\mathbb Y)\cup\{x,x'\})^4]^{1/2}\,\dd x\dd x'\nonumber\\
&\leq n\int_{\mathbb Y}\mathbb{P}(\xi_n(x,(\eta_n\cap\mathbb Y)\cup\{x\})\neq 0)^{1/2}\mathbb{E}[\xi_n(x,(\eta_n\cap\mathbb Y)\cup\{x\})^4]^{1/2}\,\dd x\nonumber\\
&\qquad
+n^2\int_{\mathbb{Y}^2}\mathbb{P}(\xi_n(x',(\eta_n\cap\mathbb Y)\cup\{x,x'\})\neq 0)^{1/4}\mathbb{P}(\xi_n(x,(\eta_n\cap\mathbb Y)\cup\{x,x'\})\neq 0)^{1/4}\nonumber\\
&\qquad\qquad\qquad\qquad\qquad\qquad\qquad\qquad\qquad
\times\mathbb{E}[\xi_n(x,(\eta_n\cap\mathbb Y)\cup\{x,x'\})^4]^{1/2}\,\dd x\dd x'\nonumber\\
&\leq n\int_{\mathbb Y}\mathbb{P}(\xi_n(x,(\eta_n\cap\mathbb Y)\cup\{x\})\neq 0)^{1/2}\mathbb{E}[\xi_n(x,(\eta_n\cap\mathbb Y)\cup\{x\})^{4+p'}]^{2/(4+p')}\,\dd x\nonumber\\
&\qquad
+n^2\int_{\mathbb{Y}^2}\mathbb{P}(\xi_n(x',(\eta_n\cap\mathbb Y)\cup\{x,x'\})\neq 0)^{1/4}\mathbb{P}(\xi_n(x,(\eta_n\cap\mathbb Y)\cup\{x,x'\})\neq 0)^{1/4}\nonumber\\
&\qquad\qquad\qquad\qquad\qquad\qquad\qquad\qquad\qquad
\times\mathbb{E}[\xi_n(x,(\eta_n\cap\mathbb Y)\cup\{x,x'\})^{4+p'}]^{2/(4+p')}\,\dd x\dd x'\nonumber\\
&\leq C_{p'}^{2/(4+p')}n\int_{\mathbb Y}\mathbb{P}(\xi_n(x,(\eta_n\cap\mathbb Y)\cup\{x\})\neq 0)^{1/2}\,\dd x\nonumber\\
&\qquad
+C_{p'}^{2/(4+p')}n^2\int_{\mathbb{Y}^2}\mathbb{P}(\xi_n(x',(\eta_n\cap\mathbb Y)\cup\{x,x'\})\neq 0)^{1/4}\mathbb{P}(\xi_n(x,(\eta_n\cap\mathbb Y)\cup\{x,x'\})\neq 0)^{1/4}
\,\dd x\dd x'\nonumber\\
&\leq C_K^{1/2}C_{p'}^{2/(4+p')}\left[n\int_{\mathbb Y}\exp(-c_K n^{\alpha_K/d}\mathrm{d}(x,K)^{\alpha_K}/2)\,\dd x
+\left(n\int_{\mathbb{Y}}\exp(-c_K n^{\alpha_K/d}\mathrm{d}(x,K)^{\alpha_K}/4)\,\dd x\right)^2\right]\nonumber\\
&\leq C'[I_{K,n}+(I_{K,n})^2]<\infty,\nonumber
\end{align}
for a suitable constant $C'\in\mathbb{R}_+$. It remains to check \eqref{eq:ineqexp}. \textcolor{black}{Since $\gamma\leq 2^{2\alpha+3}(4+{p''})$,} we have
\textcolor{black}{
\begin{align}
&\exp(-c_K n^{\alpha_K/d}\mathrm{d}(x,K)^{\alpha_K}/\gamma)\nonumber\\
&\qquad\qquad
\leq\exp\left(-\frac{{p''}\,c_{p'}}{2^{2\alpha+3}(4+{p''})}(n^{1/d}\mathrm{d}(x,K))^{\alpha_K}\right)\nonumber\\
&\qquad\qquad
=\bold{1}\{n^{1/d}\mathrm{d}(x,K)\leq 1\}\exp\left(-\frac{{p''}\,c_{p'}}{2^{2\alpha+3}(4+{p''})}(n^{1/d}\mathrm{d}(x,K))^{\alpha_K}\right)\nonumber\\
&\qquad\qquad\qquad\qquad\qquad
+\bold{1}\{n^{1/d}\mathrm{d}(x,K)>1\}\exp\left(-\frac{{p''}\,c_{p'}}{2^{2\alpha+3}(4+{p''})}(n^{1/d}\mathrm{d}(x,K))^{\alpha_K}\right)\nonumber\\
&\qquad\qquad
\leq\bold{1}\{n^{1/d}\mathrm{d}(x,K)\leq 1\}\exp\left(-\frac{{p''}\,c_{p'}}{2^{2\alpha+3}(4+{p''})}(n^{1/d}\mathrm{d}(x,K))^{\alpha}\right)\nonumber\\
&\qquad\qquad\qquad\qquad\qquad
\times\exp\left(\frac{{p''}\,c_{p'}}{2^{2\alpha+3}(4+{p''})}[(n^{1/d}\mathrm{d}(x,K))^{\alpha}-(n^{1/d}\mathrm{d}(x,K))^{\alpha_K}]\right)\nonumber\\
&\qquad\qquad\qquad\qquad\qquad
+\exp\left(-\frac{{p''}\,c_{p'}}{2^{2\alpha+3}(4+{p''})}(n^{1/d}\mathrm{d}(x,K))^{\alpha}\right)\nonumber\\
&\qquad\qquad
\leq\left(\bold{1}\{n^{1/d}\mathrm{d}(x,K)\leq 1\}\exp\left(\frac{{p''}\,c_{p'}}{2^{2\alpha+3}(4+{p''})}\right)+1\right)
\exp\left(-\frac{{p''}\,c_{p'}}{2^{2\alpha+3}(4+{p''})}(n^{1/d}\mathrm{d}(x,K))^{\alpha}\right)\nonumber\\
&\qquad\qquad
\leq\hat{C}\exp\left(-\frac{{p''}\,c_{p'}}{2^{2\alpha+3}(4+{p''})}(n^{1/d}\mathrm{d}(x,K))^{\alpha}\right),\nonumber
\end{align}
where we used that, for $a>1$, the function $x\mapsto a^x$ is non-decreasing on $[0,\infty)$.}
\\
\noindent{\it Step\,\,2:\,\,Checking\,\,the\,\,relations\,\,among\,\,the\,\,operators.} For any $x_1\in\mathbb{R}^d$, we have
\begin{align}
D_{x_1}H_n&=\bold{1}\{n^{-1/d}x_1\in\mathbb Y\}\Biggl(\xi_n(n^{-1/d}x_1,(\eta_n\cup\{n^{-1/d}x_1\})\cap\mathbb Y)\nonumber\\
&\qquad\qquad
+\sum_{z\in\eta_n\cap\mathbb Y}\xi_n(z,(\eta_n\cup\{n^{-1/d}x_1\})\cap\mathbb Y)-\sum_{z\in\eta_n\cap\mathbb Y}\xi_n(z,\eta_n\cap\mathbb Y)\Biggr)\nonumber\\
&=\bold{1}\{n^{-1/d}x_1\in\mathbb Y\}\mathcal{D}_{n^{-1/d}x_1}H_n.\nonumber
\end{align}
The relation \eqref{eq:der2} can be verified by a similar computation.
\\
\noindent$\square$

\subsection{Proof of Theorem \ref{thm:ASCLTstabilizing}}\label{sec:proveThmSTAB}

We apply Theorem \ref{cor:ASCLTAPP} and start noticing that the square integrability of the functionals $H_n$, $n\geq 1$,
follows by Lemma \ref{le:malliavinderivatives}. In the next steps we check all the other conditions of Theorem \ref{cor:ASCLTAPP}.
Hereafter, ${p''}\in (0,p')$ is the constant involved in the definition of $I_{K,n}$.\\
\noindent{\it Step\,\,1:\,\,Checking\,\,Conditions\,\,\eqref{eq:os1app}\,\,and\,\,\eqref{eq:os2app}.}
By \eqref{eq:der1} and Lemma \ref{le:LSY1} , we have
\begin{align}
\sup_{n\geq 1}\sup_{x\in\mathbb R^d}\mathbb{E}[|D_{x}H_n|^{4+{p''}}]\leq c({p''}).\nonumber
\end{align}
By \eqref{eq:der2}, the inequality $|a+b|^r\leq 2^{r-1}(|a|^r+|b|^r)$, $r\geq 1$, and again Lemma \ref{le:LSY1}, we have
\begin{align}
\sup_{n\geq 1}\sup_{x_1,x_2\in\mathbb{R}^d}\mathbb{E}[|D_{x_1,x_2}^2 H_n|^{4+{p''}}]&\leq\sup_{n\geq 1}\sup_{y_1,y_2\in\mathbb{Y}}
\mathbb{E}[|\mathcal{D}_{y_1,y_2}^{2}H_n|^{4+{p''}}]\nonumber\\
&\leq 2^{3+{p''}}\sup_{n\geq 1}\sup_{y_1,y_2\in\mathbb{Y}}(
\mathbb{E}[|\mathcal{D}_{y_1,y_2}^{2,+}H_n|^{4+{p''}}]+\mathbb{E}[|\mathcal{D}_{y_1}H_n|^{4+{p''}}])\nonumber\\
&\leq 2^{4+{p''}}c({p''}).\nonumber
\end{align}
So conditions \eqref{eq:os1app} and \eqref{eq:os2app} are verified with $p=q={p''}$ and $c=2^{4+{p''}}c({p''})$. To conclude the proof
it remains to check \eqref{hyp:convgammaH}, \eqref{hyp:series1convH},
\eqref{hyp:series2convH}, \eqref{hyp:series3convH} and \eqref{hyp:series4convH} where the quantities $\bold\Lambda$, $\bold\Gamma$, $\bold\Theta$
and $\bold\Gamma_2$ involve ${p''}$ in place of $p$ and $q$. This task is accomplished in the Steps 2-6.\\
\noindent{\it Step\,\,2:\,\,Two\,\,Preliminary\,\,Inequalities.} Let $\beta\in\mathbb{R}_+$ be fixed. By \eqref{eq:der2} and Lemma
\ref{le:LSY2}, we have
\begin{align}
&\int_{\mathbb R^d}\psi_x(H_n,\beta)^2\,\dd x=\int_{\mathbb{R}^d}\left(\int_{\mathbb R^d}\mathbb{P}(D_{x_1,x}^2H_n\neq 0)^\beta\,\dd x_1\right)^2\,\dd x\nonumber\\
&\qquad\qquad
=\int_{\mathbb{R}^d}\Biggl(\int_{\mathbb R^d}\mathbb{P}\Biggl(
\bold{1}\{n^{-1/d}x_1,n^{-1/d}x\in\mathbb Y\}\mathcal{D}_{n^{-1/d}x_1,n^{-1/d}x}^2H_n\neq 0\Biggr)^\beta\,\dd x_1\Biggr)^2\,\dd x\nonumber\\
&\qquad\qquad
=n\int_{\mathbb{Y}}\Biggl(n\int_{\mathbb Y}\mathbb{P}(
\mathcal{D}_{x_1,x}^2H_n \neq 0)^\beta\,\dd x_1\Biggr)^2\,\dd x\nonumber\\
%n\int_{\mathbb{Y}}\Biggl(n\int_{\mathbb Y}\mathbb{P}\Biggl(
%(\xi_n(x,(\eta_n\cap\mathbb Y)\cup\{x,x_1\})-\xi_n(x,(\eta_n\cap\mathbb Y)\cup\{x\}))\nonumber\\
%&\qquad\qquad\qquad\qquad\qquad\qquad
%+(\xi_n(x_1,(\eta_n\cap\mathbb Y)\cup\{x_1,x\})-\xi_n(x_1,(\eta_n\cap\mathbb Y)\cup\{x_1\}))\nonumber\\
%&\qquad\qquad
%+\sum_{z\in\eta_n\cap\mathbb Y}(\xi_n(z,(\eta_n\cap\mathbb Y)\cup\{x_1,x\})-\xi_n(z,(\eta_n\cap\mathbb Y)\cup\{x_1\})\nonumber\\
%&\qquad\qquad\qquad\qquad\qquad\qquad
%-\xi_n(z,(\eta_n\cap\mathbb Y)\cup\{x\})+\xi_n(z,\eta_n\cap\mathbb Y))\neq 0\Biggr)^\beta\,\dd x_1\Biggr)^2\,\dd x\nonumber\\
&\qquad\qquad
\leq\tilde{C}_\beta n\int_{\mathbb Y}\exp\left(-\frac{c_{p'}\beta}{2^{2\alpha+1}}(n^{1/d}\dd (x,K))^\alpha\right)\,\dd x\label{eq:ubpsi}
\end{align}
and
\begin{align}
%\label{eq:psi1S}
\int_{\mathbb{R}^d}\mathbb{P}(D_x H_n\neq 0)^\beta\,\dd x
&=\int_{\mathbb R^d}\bold{1}\{n^{-1/d}x\in\mathbb Y\}\mathbb{P}(\mathcal{D}_{n^{-1/d}x}H_n\neq 0)^\beta\,\dd x\nonumber\\
&=n\int_{\mathbb Y}\mathbb{P}(\mathcal{D}_y H_n\neq 0)^\beta\,\dd y\nonumber\\
%&=n\int_{\mathbb Y}\mathbb{P}\Biggl(
%\xi_n(y,(\eta_n\cap\mathbb Y)\cup\{y\})
%+\sum_{z\in\eta_n\cap\mathbb Y}(\xi_n(z,(\eta_n\cap\mathbb Y)\cup\{y\})-\xi_n(z,\eta_n\cap\mathbb Y))\neq 0\Biggr)^\beta\,\dd y\nonumber\\
&\leq\tilde{C}_\beta n\int_{\mathbb Y}\exp\left(-\frac{c_{p'}\beta}{2^{\alpha+1}}(n^{1/d}\dd (x,K))^\alpha\right)\,\dd x\nonumber\\
&\leq\tilde{C}_\beta n\int_{\mathbb Y}\exp\left(-\frac{c_{p'}\beta}{2^{2\alpha+1}}(n^{1/d}\dd (x,K))^\alpha\right)\,\dd x.\label{eq:psi1S}
\end{align}
Since  $\psi_x(G,\beta)$ is a non-increasing function of $\beta$, from relation \eqref{eq:ubpsi}
we get
\begin{align}
&\int_{\mathbb R^d}\psi_x(H_n,{p''}/(4(4+{p''})))^2\,\dd x\bigvee\int_{\mathbb R^d}\psi_x(H_n,{p''}/(2(4+{p''})))^2\,\dd x\nonumber\\
&=\int_{\mathbb R^d}\psi_x(H_n,{p''}/(4(4+{p''})))^2\,\dd x\leq\tilde{C}_{{p''}/(4(4+{p''}))}\,I_{K,n}.\label{eq:ubpsiNEW}
\end{align}
Similarly, from  \eqref{eq:psi1S}, since  $\mathbb{P}(D_x H_n\neq 0)^{\beta}$ is non-increasing in $\beta$ we obtain
%:Relation \eqref{eq:ubpsi}
\begin{align}
\int_{\mathbb{R}^d}\mathbb{P}(D_x H_n\neq 0)^{(1+{p''})/(4+{p''})}\,\dd x&\leq\int_{\mathbb{R}^d}\mathbb{P}(D_x H_n\neq 0)^{{p''}/(2(4+{p''}))}\,\dd x\nonumber\\
&\leq\int_{\mathbb{R}^d}\mathbb{P}(D_x H_n\neq 0)^{{p''}/(4(4+{p''}))}\,\dd x\leq\tilde{C}_{{p''}/(4(4+{p''}))}\,I_{K,n}.\label{eq:psi1SNEW}
\end{align}
\noindent{\it Step\,\,3:\,\,Checking\,\,Conditions\,\,\eqref{hyp:convgammaH}\,\,and\,\,\eqref{hyp:series1convH}.}
From now on in this proof, $C>0$ denotes a generic positive constant
(not depending on $n$) which may vary from line to line. Note that by the assumptions it follows
\begin{equation}\label{eq:variancebound}
\mathbb{V}ar(H_n)\geq C n^\tau,\quad\text{for all $n$ large enough}
\end{equation}
and
\begin{equation}\label{eq:IKbound}
I_{K,n}\leq C n^\tau,\quad\text{for all $n$ large enough.}
\end{equation}
%Similarly, by Lemma \ref{le:malliavinderivatives}, Lemma 5.2 in \cite{LSY},
%the inequality $(5.8)$ of Lemma 5.10 in \cite{LSY} and \eqref{eq:IKbound}
%FINIRE DI CORREGGERE
%Therefore,
By \eqref{eq:variancebound}, \eqref{eq:ubpsiNEW} and \eqref{eq:IKbound}, we have
\begin{align}
\bold{\Gamma}_1(H_n)&=\frac{1}{\mathbb{V}ar(H_n)}\sqrt{\int_{\mathbb{R}^d}\psi_x(H_n,{p''}/(4(4+{p''})))^2\,\dd x
+\int_{\mathbb{R}^d}\psi_x(H_n,{p''}/(2(4+{p''})))^2\,\dd x}\label{eq:gamma1n}\\
&\leq\frac{C}{n^\tau}\sqrt{\int_{\mathbb R^d}\psi_x(H_n,{p''}/(4(4+{p''})))^2\,\dd x+
\int_{\mathbb R^d}\psi_x(H_n,{p''}/(2(4+{p''})))^2\,\dd x}\nonumber\\
&\leq Cn^{-\tau/2},\quad\text{for all $n$ large enough.}\label{eq:gamma1totaledgeS}
\end{align}
Similarly, by \eqref{eq:variancebound}, \eqref{eq:psi1SNEW} and \eqref{eq:IKbound}, we have
\begin{align}
\bold{\Gamma}_2(H_n)&=\frac{1}{\mathbb{V}ar(H_n)^{3/2}}\int_{\mathbb{R}^d}\mathbb{P}(D_x H_n\neq 0)^{(1+{p''})/(4+{p''})}\,\dd x\label{eq:gamma2n}\\
&\leq\frac{C}{n^{3\tau/2}}\int_{\mathbb R^d}\mathbb{P}(D_x H_n\neq 0)^{(1+{p''})/(4+{p''})}\,\dd x\nonumber\\
&\leq C n^{-\tau/2},\quad\text{for all $n$ large enough.}\label{eq:gamma2totaledgeS}
\end{align}
Relation \eqref{hyp:convgammaH} follows immediately by  \eqref{eq:gamma1totaledgeS} and \eqref{eq:gamma2totaledgeS}. By
\eqref{eq:gamma1totaledgeS} and \eqref{eq:gamma2totaledgeS} we also have, for some $\bar n\geq 2$ large enough,
\begin{align}
\sum_{n\geq\bar n}\frac{1}{n(\log n)^2}\sum_{k=\bar n}^{n}\frac{\bold{\Gamma}_i(H_k)}{k}&\leq
C\sum_{n\geq\bar n}\frac{1}{n(\log n)^2}\sum_{k\geq\bar n}\frac{1}{k^{1+\tau/2}}<\infty,\quad\text{$i=1,2$}.\nonumber
\end{align}
Therefore \eqref{hyp:series1convH}  follows since $\bold{\Gamma}_i(H_n)<\infty$, for any $i=1,2$ and $n\geq 1$
as a consequence of the following assumptions:  $\mathbb{V}ar(H_n)\in\mathbb R_+$ and $I_{K,n}\in\mathbb R_+$ for  $n\geq 1$, \eqref{eq:gamma1n}, \eqref{eq:gamma2n},
\eqref{eq:ubpsiNEW} and \eqref{eq:psi1SNEW} (which hold for any $n\geq 1$).\\
\noindent{\it Step\,\,4:\,\,Checking\,\,Condition\,\,\eqref{hyp:series2convH}.}
By \eqref{eq:variancebound}, \eqref{eq:psi1SNEW} and \eqref{eq:IKbound}, we have
\begin{align}
\bold{\Theta}(H_k,H_l)&=\frac{1}{\sqrt{\mathbb{V}ar(H_k)\mathbb{V}ar(H_l)}}
\int_{\mathbb{R}^d}\mathbb{P}(D_x H_k\neq 0)^{{p''}/(4(4+{p''}))}\mathbb{P}(D_x H_l\neq 0)^{{p''}/(4(4+{p''}))}\,\mathrm{d}x\nonumber\\
&\leq\frac{1}{\sqrt{\mathbb{V}ar(H_k)\mathbb{V}ar(H_l)}}
\int_{\mathbb{R}^d}\mathbb{P}(D_x H_k\neq 0)^{{p''}/(4(4+{p''}))}\,\mathrm{d}x\label{eq:tetakl}\\
&\leq\frac{C}{\sqrt{k^\tau}\sqrt{l^\tau}}
\int_{\mathbb{R}^d}\mathbb{P}(D_x H_k\neq 0)^{{p''}/(4(4+{p''}))}\,\mathrm{d}x\nonumber\\
&\leq C\left(\frac{k}{l}\right)^{\tau/2},\quad\text{for all $k,l$ large enough.}\label{eq:ineqintder1}
\end{align}
Therefore, for some $\bar k\geq 2$ large enough,
\begin{align}
\sum_{n\geq\bar k}\frac{1}{n(\log n)^3}\sum_{l=\bar k}^{n}\sum_{k=\bar k}^{l}\frac{\bold{\Theta}(H_k,H_l)}{kl}
&\leq C\sum_{n\geq\bar k}\frac{1}{n(\log n)^3}\sum_{l=\bar k}^{n}\frac{1}{l^{1+\tau/2}}\sum_{k=\bar k}^{l}k^{\tau/2-1}\nonumber\\
&\leq C\sum_{n\geq\bar k}\frac{1}{n(\log n)^3}\sum_{l=\bar k}^{n}\frac{1}{l^{1+\tau/2}}
\sum_{k=2}^{l}\int_{k-1}^{k}k^{\tau/2-1}\,\dd x\nonumber\\
&\leq C\sum_{n\geq\bar k}\frac{1}{n(\log n)^3}\sum_{l=\bar k}^{n}\frac{1}{l^{1+\tau/2}}
\int_{1}^{l}x^{\tau/2-1}\,\dd x\label{eq:crucialless1}\\
&=C\sum_{n\geq\bar k}\frac{1}{n(\log n)^3}\sum_{l=\bar k}^{n}\frac{1}{l^{1+\tau/2}}
(l^{\tau/2}-1)\nonumber\\
&\leq C\sum_{n\geq\bar k}\frac{1}{n(\log n)^3}\sum_{l=\bar k}^{n}\frac{1}{l}\nonumber\\
&\leq C\sum_{n\geq\bar k}\frac{1}{n(\log n)^2}<\infty,\label{eq:suminf1}
\end{align}
where in \eqref{eq:crucialless1} we used that $0<\tau<2$. By \eqref{eq:variancebound}, the Cauchy-Schwarz inequality, \eqref{eq:psi1SNEW}
and \eqref{eq:IKbound}, for $\bar k\geq 2$ large enough and $n\geq\bar k$,
\begin{align}
&\sum_{l=\bar k}^{n}\sum_{k=1}^{\bar k-1}\frac{\bold{\Theta}(H_k,H_l)}{kl}\nonumber\\
&\quad
\leq C\sum_{l=\bar k}^{n}\sum_{k=1}^{\bar k-1}\frac{1}{k\mathbb{V}ar(H_k)^{1/2}l^{1+\tau/2}}
\left(\int_{\mathbb{R}^d}\mathbb{P}(D_x H_k\neq 0)^{{p''}/(2(4+{p''}))}\,\dd x\right)^{1/2}
\left(\int_{\mathbb R^d}\mathbb{P}(D_x H_l\neq 0)^{{p''}/(2(4+{p''}))}\,\mathrm{d}x\right)^{1/2}\nonumber\\
&\quad
\leq C\sum_{k=1}^{\bar k-1}\frac{\sqrt{I_{K,k}}}{k\mathbb{V}ar(H_k)^{1/2}}\sum_{l=\bar k}^{n}\frac{1}{l}\nonumber\\
%&\quad
%\leq C\sum_{k=1}^{\bar k-1}\frac{1}{k^{1-\tau/2}\mathbb{V}ar(H_k)^{1/2}}\sum_{l=\bar k}^{n}\frac{1}{l}\nonumber\\
&\quad
\leq C\log n.\label{eq:sum2}
\end{align}
Relations \eqref{eq:suminf1} and \eqref{eq:sum2} yield, for $\bar k\geq 2$ large enough,
\[
\sum_{n\geq\bar k}\frac{1}{n(\log n)^3}\sum_{l=\bar k}^{n}\sum_{k=1}^{l}\frac{\bold{\Theta}(H_k,H_l)}{kl}<\infty.
\]
This relation implies \eqref{hyp:series2convH} since $\bold{\Theta}(H_k,H_l)<\infty$ for any $k,l\geq 1$
as a consequence of the following assumptions: $\mathbb{V}ar(H_n)\in\mathbb R_+$ and $I_{K,n}\in\mathbb R_+$ for any $n\geq 1$, \eqref{eq:tetakl} and
\eqref{eq:psi1SNEW} (which hold for any $k,l\geq 1$).\\
\noindent{\it Step\,\,5:\,\,Checking\,\,Condition\,\,\eqref{hyp:series3convH}.} For $k,l,m,n\geq 1$, using that $\psi_x(G,\beta)$ is non-increasing in $\beta$
and the Cauchy-Schwarz inequality, we have
\begin{align}
%&\bold{\Lambda}(H_k,H_l,H_m,H_n)^2\leq
&\bold{\Gamma}(H_k,H_l,H_m,H_n)^2\nonumber\\
&\,\,\,\,\,\leq\frac{1}
{\sqrt{\mathbb{V}ar(H_k)\mathbb{V}ar(H_l)\mathbb{V}ar(H_m)\mathbb{V}ar(H_n)}}
\int_{\mathbb{R}^d}\psi_x(H_k,{p''}/(4(4+{p''})))\psi_x(H_l,{p''}/(4(4+{p''})))\,\dd x\nonumber\\
&\qquad+\frac{1}{\mathbb{V}ar(H_k)\mathbb{V}ar(H_l)}
\int_{\mathbb{R}^d}\psi_x(H_k,{p''}/(2(4+{p''})))\psi_x(H_l,{p''}/(2(4+{p''})))\,\dd x\nonumber\\
&\,\,\,\,\,\leq\left(\frac{1}
{\sqrt{\mathbb{V}ar(H_k)\mathbb{V}ar(H_l)\mathbb{V}ar(H_m)\mathbb{V}ar(H_n)}}+\frac{1}{\mathbb{V}ar(H_k)\mathbb{V}ar(H_l)}\right)\nonumber\\
&\qquad\qquad\qquad\qquad\times
\int_{\mathbb{R}^d}\psi_x(H_k,{p''}/(4(4+{p''})))\psi_x(H_l,{p''}/(4(4+{p''})))\,\dd x\nonumber\\
&\,\,\,\,\,\leq\left(\frac{1}
{\sqrt{\mathbb{V}ar(H_k)\mathbb{V}ar(H_l)\mathbb{V}ar(H_m)\mathbb{V}ar(H_n)}}+\frac{1}{\mathbb{V}ar(H_k)\mathbb{V}ar(H_l)}\right)\nonumber\\
&\qquad\qquad\qquad\qquad\times
\left(\int_{\mathbb{R}^d}\psi_x(H_k,{p''}/(4(4+{p''})))^2\,\dd x\right)^{1/2}
\left(\int_{\mathbb{R}^d}\psi_x(H_l,{p''}/(4(4+{p''})))^2\,\dd x\right)^{1/2}.\label{eq:9Dic}
\end{align}
By this relation, \eqref{eq:variancebound}, \eqref{eq:ubpsiNEW} and \eqref{eq:IKbound}, for $l,k$ large enough
\begin{align}
\bold{\Gamma}(H_k,H_l,H_k,H_l)\leq\frac{C}{(kl)^{\tau/4}}.\nonumber
\end{align}
Therefore, for some $\bar k\geq 2$ large enough
\begin{equation}\label{eq:step6}
\sum_{n\geq\bar k}\frac{1}{n(\log n)^3}\sum_{l=\bar k}^{n}\sum_{k=\bar k}^{l}\frac{\bold{\Gamma}(H_k,H_l,H_k,H_l)}{kl}
\leq C\sum_{n\geq\bar k}\frac{1}{n(\log n)^3}\sum_{l=\bar k}^{n}\frac{1}{l^{1+\tau/4}}\sum_{k=\bar k}^{l}\frac{1}{k^{1+\tau/4}}<\infty.
\end{equation}
For $\bar k\geq 2$ large enough, by \eqref{eq:9Dic}, \eqref{eq:ubpsiNEW}, \eqref{eq:IKbound} and \eqref{eq:variancebound}, we have
\begin{align}
&\sum_{n\geq\bar k}\frac{1}{n(\log n)^3}\sum_{l=\bar k}^{n}\sum_{k=1}^{\bar k-1}\frac{\bold{\Gamma}(H_k,H_l,H_k,H_l)}{kl}\nonumber\\
%&\leq C\sum_{n\geq\bar k}\frac{1}{n(\log n)^3}\sum_{l=\bar k}^{n}\sum_{k=1}^{\bar k-1}
%\frac{1}{kl}\sqrt{\frac{1}
%{\mathbb{V}ar(H_k)\mathbb{V}ar(H_l)}}
%\left(\int_{\mathbb{X}}\psi_x(H_k,{p''}/(4(4+{p''})))^2\,\dd x\right)^{1/4}
%\left(\int_{\mathbb X}\psi_x(H_l,{p''}/(4(4+{p''})))^2\,\dd x\right)^{1/4}\nonumber\\
&\leq C\sum_{n\geq\bar k}\frac{1}{n(\log n)^3}\sum_{l=\bar k}^{n}\sum_{k=1}^{\bar k-1}
\frac{1}{kl}\sqrt{\frac{I_{K,k}^{1/2}I_{K,l}^{1/2}}
	{\mathbb{V}ar(H_k)\mathbb{V}ar(H_l)}}\nonumber\\
&\leq C \left(
\sum_{k=1}^{\bar k-1}
\frac{1}{k}\sqrt{\frac{I_{K,k}^{1/2}}
	{\mathbb{V}ar(H_k)}}\right)  \sum_{n\geq\bar k}\frac{1}{n(\log n)^3}\left(\sum_{l=\bar k}^{n}
 \frac{1}{l}\sqrt{\frac{I_{K,l}^{1/2}}
	{\mathbb{V}ar(H_l)}}\right)
\nonumber\\
&\leq C\sum_{n\geq\bar k}\frac{1}{n(\log n)^3}\sum_{l=\bar k}^{n}
\frac{1}{l}\sqrt{\frac{I_{K,l}^{1/2}}
	{\mathbb{V}ar(H_l)}}\nonumber\\
&\leq C\sum_{n\geq\bar k}\frac{1}{n(\log n)^3}\sum_{l=\bar k}^{n}
\frac{1}{l^{1+\tau/4}}<\infty.\nonumber
\end{align}
This relation, along with \eqref{eq:step6}, yields \eqref{hyp:series3convH} since $\bold{\Gamma}(H_k,H_l,H_k,H_l)<\infty$ for any $k,l\geq 1$
as a consequence of the following assumptions: $\mathbb{V}ar(H_n)\in\mathbb R_+$ and $I_{K,n}\in\mathbb R_+$ for any  $n\geq 1$,  \eqref{eq:9Dic} and
\eqref{eq:ubpsiNEW} (which hold for any $k,l\geq 1$).\\
\noindent{\it Step\,\,6:\,\,Checking\,\,Condition\,\,\eqref{hyp:series4convH}.}
\textcolor{black}{By the trivial inequality
\[
\bold{\Lambda}(H_k,H_l,H_m,H_n)\leq\bold{\Gamma}(H_k,H_l,H_m,H_n),
\]}
\eqref{eq:9Dic}, the inequality $(a+b)^{1/2}\leq a^{1/2}+b^{1/2}$, $a,b\geq 0$,
and \eqref{eq:ubpsiNEW}, for $k,l\geq 1$, we have
\begin{align}
\bold{\Lambda}_{max}(H_k,H_l)&\leq C\left(\frac{1}{\mathbb{V}ar(H_k)^{3/4}\mathbb{V}ar(H_l)^{1/4}}+\frac{1}{\mathbb{V}ar(H_k)}\right)I_{K,k}^{1/2}\nonumber\\
&\qquad
\vee\left(\frac{1}{\mathbb{V}ar(H_l)^{3/4}\mathbb{V}ar(H_k)^{1/4}}+\frac{1}{\mathbb{V}ar(H_l)}\right)I_{K,l}^{1/2}\nonumber\\
&\qquad
\vee\left(\frac{1}{(\mathbb{V}ar(H_k)\mathbb{V}ar(H_l))^{1/2}}+\frac{1}{\mathbb{V}ar(H_k)}\right)I_{K,k}^{1/2}\nonumber\\
&\qquad
\vee\left(\frac{1}{(\mathbb{V}ar(H_k)\mathbb{V}ar(H_l))^{1/2}}+\frac{1}{\mathbb{V}ar(H_l)}\right)I_{K,l}^{1/2}\nonumber\\
&\qquad
\vee\left(\frac{1}{\mathbb{V}ar(H_k)^{3/4}\mathbb{V}ar(H_l)^{1/4}}+\frac{1}{(\mathbb{V}ar(H_k)\mathbb{V}ar(H_l))^{1/2}}\right)(I_{K,k}I_{K,l})^{1/4}\nonumber\\
&\qquad
\vee\left(\frac{1}{\mathbb{V}ar(H_l)^{3/4}\mathbb{V}ar(H_k)^{1/4}}+\frac{1}{(\mathbb{V}ar(H_k)\mathbb{V}ar(H_l))^{1/2}}\right)(I_{K,k}I_{K,l})^{1/4}\nonumber\\
&\qquad
\vee\frac{(I_{K,k}I_{K,l})^{1/4}}{(\mathbb{V}ar(H_k)\mathbb{V}ar(H_l))^{1/2}}\nonumber\\
&\textcolor{black}{=}C\left(\frac{I_{K,k}^{1/2}}{\mathbb{V}ar(H_k)^{3/4}}\frac{1}{\mathbb{V}ar(H_l)^{1/4}}+\frac{I_{K,k}^{1/2}}{\mathbb{V}ar(H_k)}\right)\nonumber\\
&\qquad
\vee\left(\frac{1}{\mathbb{V}ar(H_k)^{1/4}}\frac{I_{K,l}^{1/2}}{\mathbb{V}ar(H_l)^{3/4}}+\frac{I_{K,l}^{1/2}}{\mathbb{V}ar(H_l)}\right)\nonumber\\
&\qquad
\vee\left(\frac{I_{K,k}^{1/2}}{(\mathbb{V}ar(H_k))^{1/2}}\frac{1}{(\mathbb{V}ar(H_l))^{1/2}}+\frac{I_{K,k}^{1/2}}{\mathbb{V}ar(H_k)}\right)\nonumber\\
&\qquad
\vee\left(\frac{1}{(\mathbb{V}ar(H_k))^{1/2}}\frac{I_{K,l}^{1/2}}{(\mathbb{V}ar(H_l))^{1/2}}+\frac{I_{K,l}^{1/2}}{\mathbb{V}ar(H_l)}\right)\nonumber\\
&\qquad
\vee\left(\frac{(I_{K,k}I_{K,l})^{1/4}}{\mathbb{V}ar(H_k)^{3/4}\mathbb{V}ar(H_l)^{1/4}}+\frac{(I_{K,k}I_{K,l})^{1/4}}{(\mathbb{V}ar(H_k)\mathbb{V}ar(H_l))^{1/2}}\right)\nonumber\\
&\qquad
\vee\left(\frac{(I_{K,k}I_{K,l})^{1/4}}{\mathbb{V}ar(H_l)^{3/4}\mathbb{V}ar(H_k)^{1/4}}+\frac{(I_{K,k}I_{K,l})^{1/4}}{(\mathbb{V}ar(H_k)\mathbb{V}ar(H_l))^{1/2}}\right)\nonumber\\
&\qquad
\vee\frac{(I_{K,k}I_{K,l})^{1/4}}{(\mathbb{V}ar(H_k)\mathbb{V}ar(H_l))^{1/2}}.\label{eq:10Dicdp}
\end{align}
Consequently, by \eqref{eq:IKbound} and \eqref{eq:variancebound}, for some $\bar k\geq 2$ large enough and $l\geq k\geq\bar k$, we have
\begin{align}
\bold{\Lambda}_{max}(H_k,H_l)\leq\frac{C}{k^{\tau/2}}.\nonumber
\end{align}
Therefore, for some $\bar k\geq 2$ large enough, we have
\begin{equation}\label{eq:10Dic}
\sum_{n\geq\bar k}\frac{1}{n(\log n)^3}\sum_{l=\bar k}^{n}\sum_{k=\bar k}^{l}\frac{\bold{\Lambda}_{max}(H_k,H_l)}{kl}
\leq
C\sum_{n\geq\bar k}\frac{1}{n(\log n)^3}\sum_{l=\bar k}^{n}\frac{1}{l}\sum_{k=\bar k}^{l}\frac{1}{k^{1+\tau/2}}<\infty.
\end{equation}
Let $\bar k\geq 2$ large enough and $l\geq\bar k>k$, let $C_i(k)\in\mathbb{R}_+$, $i=0,1,\ldots,12$, denote suitable constants, depending on $k$ but not on $l$,
which may vary from line to line. By \eqref{eq:10Dicdp}, \eqref{eq:IKbound} and \eqref{eq:variancebound}, we have
\begin{align}
\bold{\Lambda}_{max}(H_k,H_l)&\leq\left(C_1(k)\frac{1}{\mathbb{V}ar(H_l)^{1/4}}+C_2(k)\right)\nonumber\\
&\vee\left(C_3(k)\frac{I_{K,l}^{1/2}}{\mathbb{V}ar(H_l)^{3/4}}+C_4(k)\frac{1}{\mathbb{V}ar(H_l)}\right)\vee\left(C_5(k)\frac{1}{(\mathbb{V}ar(H_l))^{1/2}}+C_6(k)\right)\nonumber\\
&\vee\left(C_7(k)\frac{I_{K,l}^{1/2}}{(\mathbb{V}ar(H_l))^{1/2}}+\frac{I_{K,l}^{1/2}}{\mathbb{V}ar(H_l)}\right)
\vee\left(C_8(k)\frac{(I_{K,l})^{1/4}}{\mathbb{V}ar(H_l)^{1/4}}+C_9(k)\frac{(I_{K,l})^{1/4}}{(\mathbb{V}ar(H_l))^{1/2}}\right)\nonumber\\
&\vee\left(C_{10}(k)\frac{(I_{K,l})^{1/4}}{\mathbb{V}ar(H_l)^{3/4}}+C_{11}(k)\frac{(I_{K,l})^{1/4}}{(\mathbb{V}ar(H_l))^{1/2}}\right)
\vee C_{12}(k)\frac{(I_{K,l})^{1/4}}{(\mathbb{V}ar(H_l))^{1/2}}\nonumber\\
&\leq\left(C_1(k)l^{-\tau/4}+C_2(k)\right)\vee\left(C_3(k)l^{-\tau/4}+C_4(k)l^{-\tau/2}\right)\nonumber\\
&\vee\left(C_5(k)l^{-\tau/2}+C_6(k)\right)\vee\left(C_7(k)+l^{-\tau/2}\right)\nonumber\\
&\vee\left(C_8(k)+C_9(k)l^{-\tau/4}\right)\vee\left(C_{10}(k)l^{-\tau/2}+C_{11}(k)l^{-\tau/4}\right)\vee C_{12}(k)l^{-\tau/4}\nonumber\\
&\leq C_0(k).\nonumber
\end{align}
Therefore, for $\bar k\geq 2$ large enough, we have
\begin{align*}
\sum_{n\geq\bar k}\frac{1}{n(\log n)^3}\sum_{l=\bar k}^{n}\sum_{k=1}^{\bar k-1}\frac{\bold{\Lambda}_{max}(H_k,H_l)}{kl}
&\leq C_0(k)\sum_{n\geq\bar k}\frac{1}{n(\log n)^3}\sum_{l=\bar k}^{n}\frac{1}{l}<\infty.
\end{align*}
This relation, along with  \eqref{eq:10Dic}, yields \eqref{hyp:series4convH} since $\bold{\Lambda}_{max}(H_k,H_l)<\infty$ for any $k,l\geq 1$
as a consequence of the following assumptions: $\mathbb{V}ar(H_n)\in\mathbb R_+$ and $I_{K,n}\in\mathbb R_+$ for any $n\geq 1$, and \eqref{eq:10Dicdp} (which holds for any $k,l\geq 1$).
\textcolor{black}{\section{Proofs of Corollaries \ref{cor:totaledge}, \ref{cor:cliquecount} and \ref{cor:voronoimosaic}}\label{sec:rg}}

\subsection{Proof of Corollary \ref{cor:totaledge}}

Set $H_n:=n^{m/d}L_{NG_k}^{(m)}(\eta_n\cap\mathbb Y)$ and note that, for $F_n$ defined in
\eqref{eq:totaledge}, we have
\[
F_n=\frac{H_n-\mathbb{E}[H_n]}{\sqrt{\mathbb{V}ar(H_n)}}.
\]
The claim follows by Theorem \ref{thm:ASCLTstabilizing} with $\tau=1$. Indeed: $(i)$ By Theorem 2.1 and Lemma 6.3 in \cite{PY} (see also the discussion
on the related variance bounds in \cite{LSY}) we have $n^{1-2m/d}=O(\mathbb{V}ar(L_{NG_k}^{(m)}(\eta_n\cap\mathbb Y))$,
%$\mathbb{V}ar(L_{NG_k}^{(q)}(\eta_n\cap\mathbb Y))=\Theta(n^{1-2q/d})$
and so $n=O(\mathbb{V}ar(H_n))$. $(ii)$ One can show that the functionals $H_n$, $n\geq 1$, are stabilizing following the proof of Theorem 3.1 in \cite{LSY}.
We briefly sketch the line of the proof. For $n\geq 1$ and $x\in\eta_n\cap\mathbb Y$, we set $\xi_n(x,\eta_n\cap\mathbb Y):=n^{m/d}\xi^{(m)}(x,\eta_n\cap\mathbb Y)$.
It turns out that $H_n=\sum_{x\in\eta_n\cap\mathbb Y}\xi_n(x,\eta_n\cap\mathbb Y)$ and that
the score functions $\xi_n$, $n\geq 1$, satisfy \eqref{eq:expstab}, \eqref{eq:boundedmoment} and \eqref{eq:expdecay} with $K:=\mathbb Y$.
%(see the proof of Theorem 3.1 in \cite{LSY}).
$(iii)$ By \eqref{eq:IK} easily follows $I_{K,n}=I_{\mathbb Y,n}=n\ell_d(\mathbb Y)$.

\subsection{Proof of Corollary \ref{cor:cliquecount}}

Setting $H_n:=\mathcal{C}_k(\eta_n\cap\mathbb Y,rn^{-1/d})$, $n\geq 1$, the claim follows
by Theorem \ref{thm:ASCLTstabilizing} with $\tau=1$. Indeed: $(i)$ As noticed in the proof of Theorem 2.5 in Section \textcolor{black}{7} of
\cite{PY2} (see also Remark$(iv)$ on p. \textcolor{black}{967} in \cite{LSY}), we have $\inf_{n\geq 1}\mathbb{V}ar(H_n)/n>0$,
and so $n=O(\mathbb{V}ar(H_n))$. $(ii)$ One can show that the functionals $H_n$, $n\geq 1$, are stabilizing following the proof of Theorem 3.15 in \cite{LSY}.
We briefly sketch the line of the proof. For $n\geq 1$ and $x\in\eta_n\cap\mathbb Y$, we denote by $\xi_n(x,\eta_n\cap\mathbb Y)$
the number of cliques of order $k+1$ in $G(\eta_n\cap\mathbb Y,rn^{-1/d})$ containing $x$.
It turns out that $H_n=\sum_{x\in\eta_n\cap\mathbb Y}\xi_n(x,\eta_n\cap\mathbb Y)$ and that
the score functions $\xi_n$, $n\geq 1$, satisfy \eqref{eq:expstab}, \eqref{eq:boundedmoment} and \eqref{eq:expdecay} with $K:=\mathbb Y$.
%(see the proof of Theorem 3.15 in \cite{LSY}).
$(iii)$ By \eqref{eq:IK} easily follows $I_{K,n}=I_{\mathbb Y,n}=n\ell_d(\mathbb Y)$.

\subsection{Proof of Corollary \ref{cor:voronoimosaic}}

Set $H_n:=n(\ell_d(A_n)-\ell_d(A))$, $n\geq 1$, and note that, for $F_n$ defined in \eqref{eq:voronoi},
%with $\varphi(\cdot):=\ell_d(\cdot)$,
we have
\[
F_n=\frac{H_n-\mathbb{E}[H_n]}{\sqrt{\mathbb{V}ar(H_n)}}.
\]
The claim follows by applying Theorem \ref{thm:ASCLTstabilizing} with $\tau:=1-1/d$. Indeed:
$(i)$ By Theorem 1.2 in \cite{S} (see also Theorem 1.1 in \cite{TY})
%and Theorem 2.2 in \cite{Y}
we have $n^{1-1/d}=O(\mathbb{V}ar(H_n))$. $(ii)$ One can show that the functionals $H_n$, $n\geq 1$, are stabilizing following
the proof of Theorem 3.4 in \cite{LSY} (and the references cited therein). For the sake of completeness, we briefly sketch the line of the proof. For $n\geq 1$ and $x\in\eta_n\cap\mathbb Y$,
%and
%a fixed $\mathbb Y'\subseteq\mathbb Y$,
we set
%\[
%\xi(x,\eta_n\cap\mathbb Y):=\bold{1}_{\mathbb Y'}(x)\ell_d(\mathrm{C}(x,\eta_n\cap\mathbb Y)\cap\mathbb Y^{'c})
%-\bold{1}_{\mathbb Y^{'c}}(x)\ell_d(\mathrm{C}(x,\eta_n\cap\mathbb Y)\cap\mathbb Y')
%\]
\[
\xi(x,\eta_n\cap\mathbb Y):=\bold{1}_{A}(x)\ell_d(\mathrm{C}(x,\eta_n\cap\mathbb Y)\cap A^{c})
-\bold{1}_{A^{c}}(x)\ell_d(\mathrm{C}(x,\eta_n\cap\mathbb Y)\cap A)
\]
and $\xi_n(x,\eta_n\cap\mathbb Y):=n\xi(x,\eta_n\cap\mathbb Y)$. It turns out that $H_n=\sum_{x\in\eta_n\cap\mathbb Y}\xi_n(x,\eta_n\cap\mathbb Y)$ and that
the score functions $\xi_n$, $n\geq 1$, satisfy \eqref{eq:expstab}, \eqref{eq:boundedmoment} and \eqref{eq:expdecay} with $K:=\partial A$.
$(iii)$ To prove $I_{K,n}=I_{\partial A,n}=O(n^{1-1/d})$, for the sake of clarity, we reproduce the argument
in the proof of Theorem 2.3 in \cite{LSY}. For all the notions and results of geometric measure theory considered hereon we refer the reader to
\cite{AFP}. We preliminary note that since $A$ is a convex body, it has a $(d-1)$-rectifiable boundary (i.e., its boundary is the
Lipschitz image of a bounded set in $\mathbb{R}^{d-1}$). Then the $(d-1)$-dimensional upper Minkowski content of $\partial A$, denoted by $\overline{\mathcal{M}}^{d-1}(\partial A)$, is
a scalar multiple of the $(d-1)$-dimensional Hausdorff measure of $\partial A$, denoted by $\mathcal{H}^{d-1}(\partial A)$. Therefore
$\overline{\mathcal{M}}^{d-1}(\partial A)<\infty$, and by Lemma 5.12 in \cite{LSY}, there exists a positive constant $C>0$ such that
\begin{equation}\label{eq:Hausdorffineq}
\mathcal{H}^{d-1}(\partial A_r)\leq C(1+r^{d-1}),\quad r>0
\end{equation}
where $\partial A_r:=\{x\in\mathbb{R}^d:\,\,\mathrm{d}(x,\partial A)\leq r\}$. Setting $c:={p''}\,c_{p'}/[2^{2\alpha+3}(4+{p''})]$, we finally have
\textcolor{black}{
\begin{align}
I_{K,n}&=n\int_{\mathbb Y}\exp(-c n^{\alpha/d}\mathrm{d}(x,\partial A)^\alpha)\,\dd x\nonumber\\
%&=n\int_{\partial A}\exp(-c n^{\alpha/d}\mathrm{d}(x,\partial A)^\alpha)\,\dd x+
&=n\int_{\mathbb Y\setminus\partial A}\exp(-c n^{\alpha/d}\mathrm{d}(x,\partial A)^\alpha)\,\dd x\nonumber\\
%&\leq n\ell_d(\partial A)+n\int_{\mathbb Y\setminus\partial A}\exp(-c n^{\alpha/d}\mathrm{d}(x,\partial A)^\alpha)\,\dd x\nonumber\\
%&n\int_{\mathbb Y\setminus\partial A}\exp(-c n^{\alpha/d}\mathrm{d}(x,\partial A)^\alpha)\,\dd x\nonumber\\
&=n\int_0^{\infty}\exp(-c n^{\alpha/d}r^\alpha)\mathcal{H}^{d-1}(\partial A_r)\,\dd r\label{eq:coarea}\\
&\leq C n\int_0^{\infty}\exp(-c n^{\alpha/d}r^\alpha)(1+r^{d-1})\,\dd r\label{eq:lemma5.12}\\
&\leq C n^{1-1/d}\int_0^\infty\mathrm{e}^{-c u^\alpha}(1+u^{d-1})\,\mathrm{d}u,\nonumber
\end{align}
}
where \eqref{eq:coarea} follows by the coarea formula and \eqref{eq:lemma5.12} is a consequence of \eqref{eq:Hausdorffineq}.

\section{Proof of Theorem \ref{cor:ASCLTAPP}}\label{sec:poissonspace}

The proof of Theorem \ref{cor:ASCLTAPP} relies on some  inequalities for Malliavin's operators on the Poisson space recently derived by Last, Peccati and Schulte \cite{LPS}. We also exploit a general characterization of the almost sure version of the classical weak convergence of random variables due to Ibragimov and Lifshits \cite{IL}.
In Subsection \ref{sec:prelMall} we provide some further preliminaries about the Malliavin calculus on the Poisson space, in Subsection \ref{sec:prellemma}
we give the preliminary lemmas which will be exploited to prove Theorem \ref{cor:ASCLTAPP} and finally in Subsection \ref{sec:proofthm} we prove the theorem.

\subsection{Elements of the Malliavin calculus on the Poisson space}\label{sec:prelMall}

For $r\in\mathbb{R}_+$ and $n\geq 1$ integer, we denote by $L^r(\mu^{\otimes n})$ the set of all measurable functions $g:\mathbb X^n\to\mathbb R$ such that
$\int_{\mathbb X^n}|g(x)|^r\mu^{\otimes n}(\dd x)<\infty$. We call a function $g:\mathbb X^n\to\mathbb R$ symmetric if it is invariant under permutations of its arguments, and denote
by $L_s^2(\mu^{\otimes n})$ the set of all symmetric functions  $g\in L^2(\mu^{\otimes n})$. For $g_1,g_2\in L^2(\mu^{\otimes n})$, we define $\langle g_1,g_2\rangle_n:=\int_{\mathbb X^n}g_1(x)g_2(x)\mu^{\otimes n}(\dd x)$
and $\|g_1\|_n:=\sqrt{\langle g_1,g_1\rangle_n}$.

For $F\in L_{\eta}^2$, $F=f(\eta)$, we extend the definition of the Malliavin operators $D$ and $D^2$ defining
$D_{x_1,\ldots,x_n}^n F:=D_{x_n}(D_{x_1,\ldots,x_{n-1}}^{n-1}F)$, $x_1,\ldots,x_n\in\mathbb X$, $n\geq 3$.
It is well-known that every $F\in L_{\eta}^2$ admits the representation
\begin{equation}\label{eq:chaos}
F=\EE[F]+\sum_{n\geq 1}I_n(g_n),
\end{equation}
where $g_n(x_1,\ldots,x_n)=\frac{1}{n!}\EE[D_{x_1,\ldots,x_n}^n F]$ and,
for $g_n\in L_s^2(\mu^{\otimes n})$, we denote by $I_n(g_n)$ the $n$th order Wiener-It$\mathrm{\hat{o}}$ integral with respect
to the centered Poisson measure $\eta(\dd x)-\mu(\dd x)$, see e.g. \cite{LastPen}. Another operator that we
shall consider is the so-called Ornstein-Uhlenbeck generator $L$. Given $F\in L_{\eta}^2$ of the form \eqref{eq:chaos},
we write $F\in\mathrm{Dom}(L)$ if $\sum_{n\geq 1}n^2 n!\|g_n\|_n^2<\infty$. In this case
we define
\[
L F:=-\sum_{n\geq 1}n I_n(g_n).
\]
The (pseudo) inverse of $L$ is given by
\[
L^{-1}F:=-\sum_{n\geq 1}n^{-1}I_n(g_n).
\]
It can be easily checked that the random variable $L^{-1}F$ is a well-defined element of $L_{\eta}^2$ for every $F\in L_{\eta}^2$.
In the following, we write $F\in\mathrm{Dom}(D)$ if $F\in L_{\eta}^2$ and $\int_{\mathbb X}\EE[|D_x F|^2]\mu(\dd x)<\infty$.

\subsection{Preliminary lemmas}\label{sec:prellemma}

In this subsection we provide some relations among Malliavin's operators on the Poisson space, which will be
crucial to prove Theorem \ref{cor:ASCLTAPP}. The following lemmas hold.

\begin{Lemma}\label{le:DLinverse}
For $F\in L_\eta^2$ and $r\geq 1$, we have
\[
\mathbb{E}[|D_x L^{-1}F|^r]\leq\mathbb{E}[|D_xF|^r],\quad\text{$\mu$-a.e. $x\in\mathbb X$.}
\]
\end{Lemma}

\begin{Lemma}\label{le:covariance}
For any $F,G\in\mathrm{Dom}(D)$ such that $\mathbb{E}[F]=0$, we have
\[
\mathbb{E}[FG]=\mathbb{E}[\langle DG,-DL^{-1}F\rangle_{1}].
\]
\end{Lemma}

\begin{Lemma}\label{le:boundvariace}
For $F\in\mathrm{Dom}(D)$ such that $\EE[F]=0$, we have
\[
\mathbb{V}\mathrm{ar}(\langle DF,-DL^{-1}F\rangle_{1})\leq\gamma_1(F)^2,
\]
where
\begin{align}
\gamma_1(F)^2&:=4\int_{\mathbb{X}^3}(\mathbb{E}[(D_{x_1,x_3}^{2}F)^2(D_{x_2,x_3}^{2}F)^2])^{1/2}(\mathbb{E}[(D_{x_1}F)^2(D_{x_2}F)^2])^{1/2}\mu^{\otimes 3}(\mathrm{d}x_1,\dd x_2,\dd x_3)\nonumber\\
&\qquad\qquad\qquad
+\int_{\mathbb{X}^3}\mathbb{E}[(D_{x_1,x_3}^{2}F)^2(D_{x_2,x_3}^{2}F)^2]\mu^{\otimes 3}(\mathrm{d}x_1,\dd x_2,\dd x_3).\label{eq:gamma1}
%\label{eq:2le}
\end{align}
\end{Lemma}
We refer the reader to \cite{LPS} for the proof of these lemmas (see, respectively, Lemma 3.4,
the first displayed formula in the proof of Proposition 4.1 and the statement of Proposition 4.1 itself, therein).

The next lemma gives a Gaussian bound for functionals of the Poisson measure. We refer the reader to Lemma 2.2 in \cite{BNT} and Lemma 3.1 in \cite{Zh}
for analogous inequalities on the Wiener and Rademacher spaces, respectively.

\begin{Lemma}\label{le:ledisfcar}
For $F\in\mathrm{Dom}(D)$ such that $\EE[F]=0$, we have
\begin{align}
|\EE[\mathrm{e}^{\bold{i}tF}]-\mathrm{e}^{-t^2/2}|\leq |t|^2(|1-\mathbb{E}[F^2]|+\gamma_1(F))+\frac{|t|^3}{\sqrt 2}\gamma_2(F),\quad\text{$\forall$ $t\in\mathbb R$}\nonumber
\end{align}
where $\bold{i}:=\sqrt{-1}$ and
\begin{equation}\label{eq:gamma2}
\gamma_2(F):=\int_{\mathbb{X}}\mathbb{E}[|D_xF|^3]\mu(\mathrm{d}x).
\end{equation}
\end{Lemma}
\noindent${\it Proof.}$ Let $\mathcal{C}_b^{1}(\mathbb R)$ be the family of bounded and differentiable functions from $\mathbb R$ to $\mathbb R$
with a bounded first derivative. For any $g\in\mathcal{C}_b^1(\mathbb R)$ it holds $g(F)\in\mathrm{Dom}(D)$ (indeed,
we clearly have $g(F)\in L_\eta^2$ and, by the mean value theorem, $|D_xg(F)|\leq\|g'\|_\infty |D_x F|$, for any $x\in\mathbb X$).
Therefore, by Lemma \ref{le:covariance}
\[
\mathbb{E}[Fg(F)]=\mathbb{E}[\langle Dg(F),-DL^{-1}F\rangle_{1}],\quad\text{$\forall$ $g\in\mathcal{C}_b^1(\mathbb R)$.}
\]
So, for any $t\in\mathbb R$,
\begin{align}
\mathbb{E}[F\mathrm{e}^{\bold{i}tF}]&=\mathbb{E}[F\cos(tF)]+\bold{i}\mathbb{E}[F\sin(tF)]\nonumber\\
&=\mathbb{E}[\langle D\cos(tF),-DL^{-1}F\rangle_{1}]+\bold{i}\mathbb{E}[\langle D\sin(tF),-DL^{-1}F\rangle_{1}].\label{eq:GB1}
\end{align}
For $F=f(\eta)$, by Taylor's formula with integral remainder, for $(x,t)\in\mathbb X\times\mathbb R$,
\begin{equation}\label{eq:GB2}
D_x\cos(tF)=\cos(tf(\eta+\varepsilon_x))-\cos(t f(\eta))=-t\sin(tF)D_xF+t^2R_{1,t}(x)
\end{equation}
and similarly
\begin{equation}\label{eq:GB3}
D_x\sin(tF)=\sin(tf(\eta+\varepsilon_x))-\sin(t f(\eta))=t\cos(tF)D_xF+t^2R_{2,t}(x),
\end{equation}
where
\[
R_{1,t}(x):=\int_{f(\eta)}^{f(\eta+\varepsilon_x)}(u-f(\eta+\varepsilon_x))\cos tu\,\mathrm{d}u
\]
and
\[
R_{2,t}(x):=\int_{f(\eta)}^{f(\eta+\varepsilon_x)}(u-f(\eta+\varepsilon_x))\sin tu\,\mathrm{d}u.
\]
Note that, for $i=1,2$, $|R_{i,t}(x)|\leq\frac{|D_xF|^2}{2}$ and so
\begin{equation}\label{eq:modulusR}
|R_{1,t}(x)+\bold{i}R_{2,t}(x)|\leq\frac{|D_x F|^2}{\sqrt 2}.
\end{equation}
By \eqref{eq:GB1}, \eqref{eq:GB2} and \eqref{eq:GB3}, we have
\begin{align}
\mathbb{E}[F\mathrm{e}^{\bold{i}tF}]&=t\mathbb{E}[(-\sin(tF)+\bold{i}\cos(tF))\langle DF,-DL^{-1}F\rangle_{1}]+t^2\mathbb{E}[\langle R_{1,t}+\bold{i}R_{2,t},-DL^{-1}F\rangle_{1}]\nonumber\\
&=\bold{i}t\mathbb{E}[\mathrm{e}^{\bold{i}tF}\langle DF,-DL^{-1}F\rangle_{1}]+t^2\mathbb{E}[\langle R_{1,t}+\bold{i}R_{2,t},-DL^{-1}F\rangle_{1}].\label{eq:derfcar}
\end{align}
We put $\varphi(t):=\mathrm{e}^{t^2/2}\EE[\mathrm{e}^{\bold{i}tF}]$, $t\in\mathbb R$. By the mean value theorem, \eqref{eq:derfcar} and \eqref{eq:modulusR}, we have
\begin{align}
|\varphi(t)-\varphi(0)|&\leq |t|\sup_{u\in [0,t]}|\varphi'(u)|\nonumber\\
&=|t|\sup_{u\in [0,t]}|u\mathrm{e}^{u^2/2}\EE[\mathrm{e}^{\bold{i}uF}]+\bold{i}\mathrm{e}^{u^2/2}\mathbb{E}[F\mathrm{e}^{\bold{i}uF}]|\nonumber\\
&=|t|\sup_{u\in [0,t]}|u\mathrm{e}^{u^2/2}\EE[\mathrm{e}^{\bold{i}uF}]-u\mathrm{e}^{u^2/2}\mathbb{E}[\mathrm{e}^{\bold{i}uF}\langle DF,-DL^{-1}F\rangle_{1}]\nonumber\\
&\quad\quad\quad+\bold{i}u^2\mathrm{e}^{u^2/2}\mathbb{E}[\langle R_{1,u}+\bold{i}R_{2,u},-DL^{-1}F\rangle_{1}]|\nonumber\\
&\leq|t|^2\mathrm{e}^{t^2/2}\sup_{u\in [0,t]}|\mathbb{E}[\mathrm{e}^{\bold{i}uF}(1-\langle DF,-DL^{-1}F\rangle_{1})]|\nonumber\\
&\qquad\qquad\qquad+|t|^3\mathrm{e}^{t^2/2}\sup_{u\in [0,t]}\mathbb{E}[|\langle R_{1,u}+\bold{i}R_{2,u},-DL^{-1}F\rangle_{1}|]\nonumber\\
&\leq|t|^2\mathrm{e}^{t^2/2}\mathbb{E}[|1-\langle DF,-DL^{-1}F\rangle_{1}|]\nonumber\\
&\qquad\qquad\qquad+\frac{|t|^3\mathrm{e}^{t^2/2}}{\sqrt 2}\int_{\mathbb X}\mathbb{E}[|D_x F|^2|D_xL^{-1}F|]\mu(\mathrm{d}x).\nonumber
\end{align}
So by Lemma \ref{le:covariance}, the Cauchy-Schwarz inequality and Lemma \ref{le:boundvariace}, it follows
\begin{align}
|\EE[\mathrm{e}^{\bold{i}tF}]-\mathrm{e}^{-t^2/2}|&\leq |t|^2|1-\mathbb{E}[F^2]|+|t|^2\mathbb{E}[|\mathbb{E}[F^2]-\langle DF,-DL^{-1}F\rangle_{1}|]\nonumber\\
&\qquad\qquad\qquad
+\frac{|t|^3}{\sqrt 2}\int_{\mathbb X}\mathbb{E}[|D_x F|^2|DL^{-1}F|]\mu(\mathrm{d}x)\nonumber\\
&\leq |t|^2|1-\mathbb{E}[F^2]|+|t|^2\sqrt{\mathbb{V}\mathrm{ar}(\langle DF,-DL^{-1}F\rangle_{1})}\nonumber\\
&\qquad\qquad\qquad
+\frac{|t|^3}{\sqrt 2}\int_{\mathbb X}\mathbb{E}[|D_x F|^2|D_xL^{-1}F|]\mu(\mathrm{d}x)\nonumber\\
&\leq |t|^2(|1-\mathbb{E}[F^2]|+\gamma_1(F))+\frac{|t|^3}{\sqrt 2}\int_{\mathbb X}\mathbb{E}[|D_x F|^2|D_xL^{-1}F|]\mu(\mathrm{d}x).\label{eq:1le}
\end{align}
By H\"older's inequality and Lemma \ref{le:DLinverse}, we have
\begin{align}
\int_X\mathbb{E}[|D_xF|^2|D_x L^{-1}F|]\mu(\mathrm{d}x)&\leq\int_{\mathbb{X}}\mathbb{E}[|D_xF|^3]^{2/3}\mathbb{E}[|D_x L^{-1}F|^3]^{1/3}\mu(\mathrm{d}x)\nonumber\\
&\leq\int_{\mathbb{X}}\mathbb{E}[|D_xF|^3]^{2/3}\mathbb{E}[|D_x F|^3]^{1/3}\mu(\mathrm{d}x)\nonumber\\
&=\int_{\mathbb{X}}\mathbb{E}[|D_xF|^3]\mu(\mathrm{d}x).\nonumber
%\label{eq:3le}
\end{align}
The claim follows combining this latter inequality with \eqref{eq:1le}.
\\
\noindent$\square$

We conclude this subsection recalling a result, due to Ibragimov and Lifshits \cite{IL}, which provides  general sufficient conditions
for the ASCLT to hold.

\begin{Lemma}\label{thm:IL}
Let $\{X_n\}_{n\geq 1}$ be a sequence of real-valued random variables converging in distribution to a random variable $X$, and put
\begin{equation}\label{eq:deltan}
\Delta_n(t):=\frac{1}{\log n}\sum_{k=1}^{n}\frac{1}{k}\left(\mathrm{e}^{\bold{i}t X_k}-\EE[\mathrm{e}^{\bold{i}t X}]\right),\quad t\in\mathbb R.
\end{equation}
If, for all $r\in\mathbb{R}_+$,
\begin{equation*}
\sup_{|t|\leq r}\sum_{n\geq 2}\frac{\mathbb{E}[|\Delta_n(t)|^2]}{n\log n}<\infty,
\end{equation*}
then \eqref{eq:ASCLT} holds with $X$ in place of $Z$. If $X$ is distributed according to the standard normal law, then $\{X_n\}_{n\geq 1}$ satisfies the ASCLT.
\end{Lemma}

\subsection{Ancillary ASCLTs on the Poisson space and proof of Theorem \ref{cor:ASCLTAPP}}\label{sec:proofthm}

Let $F_1,F_2,F_3,F_4\in L_\eta^2$ be such that $\EE[F_i]=0$, $i=1,2,3,4$. We generalize the definition of the quantity $\gamma_1(F)^2$ setting
\begin{align}
\gamma(F_1,F_2,F_3,F_4)^2:&=4\lambda(F_1,F_2,F_3,F_4)^2+\int_{\mathbb{X}^3}\mathbb{E}[(D_{x_1,x_3}^{2}F_1)^2(D_{x_2,x_3}^{2}F_2)^2]\mu^{\otimes 3}(\mathrm{d}x_1,\dd x_2,\dd x_3),\nonumber
\end{align}
where
\begin{align}
\lambda(F_1,F_2,F_3,F_4)^2:&=
\int_{\mathbb{X}^3}(\mathbb{E}[(D_{x_1,x_3}^{2}F_1)^2(D_{x_2,x_3}^{2}F_2)^2])^{1/2}(\mathbb{E}[(D_{x_1}F_3)^2(D_{x_2}F_4)^2])^{1/2}\mu^{\otimes 3}(\mathrm{d}x_1,\dd x_2,\dd x_3).\nonumber
\end{align}
Note that $\gamma_1(F_1)^2=\gamma(F_1,F_1,F_1,F_1)^2$,
\begin{equation}\label{eq:fubini}
\lambda(F_1,F_2,F_3,F_4)^2=\lambda(F_2,F_1,F_4,F_3)^2\quad\text{and}\quad\gamma(F_1,F_2,F_3,F_4)^2=\gamma(F_2,F_1,F_4,F_3)^2.
\end{equation}
We also put
\[
\theta(F_1,F_2):=\int_{\mathbb X}(\mathbb{E}[|D_x F_1|^2])^{1/2}(\EE[|D_x F_2|^2])^{1/2}\mu(\mathrm{d}x).
\]
The following ASCLT will be proved at the end of this subsection.

\begin{Proposition}\label{thm:ASCLT}
Assume $\{F_n\}_{n\geq 1}\subset L_\eta^2$ and $\mathbb{E}[F_n]=0$, $\mathbb{E}[F_n^2]=1$, for any $n\geq 1$. Moreover, suppose
\begin{equation}\label{hyp:convgamma}
\lim_{n\to\infty}\gamma_1(F_n)=\lim_{n\to\infty}\gamma_2(F_n)=0,
\end{equation}
\begin{equation}\label{hyp:series1conv}
\sum_{n\geq 2}\frac{1}{n(\log n)^2}\sum_{k=1}^{n}\frac{\gamma_i(F_k)}{k}<\infty,\quad\text{$i=1,2$}
\end{equation}
\begin{equation}\label{hyp:series2conv}
\sum_{n\geq 2}\frac{1}{n(\log n)^3}\sum_{l=1}^{n}\sum_{k=1}^{l}\frac{\theta(F_k,F_l)}{kl}<\infty,
\end{equation}
\begin{equation}\label{hyp:series3conv}
\sum_{n\geq 2}\frac{1}{n(\log n)^3}\sum_{l=1}^{n}\sum_{k=1}^{l}\frac{\gamma(F_k,F_l,F_k,F_l)}{kl}<\infty,
\end{equation}
\begin{equation}\label{hyp:series4conv}
\sum_{n\geq 2}\frac{1}{n(\log n)^3}\sum_{l=1}^{n}\sum_{k=1}^{l}\frac{\lambda_{max}(F_k,F_l)}{kl}<\infty,
\end{equation}
where
\begin{align}
\lambda_{max}(F_k,F_l):=&\lambda(F_k,F_k,F_k,F_l)\vee\lambda(F_l,F_l,F_l,F_k)\vee
\lambda(F_k,F_k,F_l,F_l)\vee\lambda(F_l,F_l,F_k,F_k)\nonumber\\
&\vee\lambda(F_l,F_k,F_k,F_k)\vee\lambda(F_k,F_l,F_l,F_l)
\vee\lambda(F_k,F_l,F_l,F_k).\nonumber
%\label{eq:lambdamax}
\end{align}
Then $\{F_n\}_{n\geq 1}$ satisfies the ASCLT.
\end{Proposition}

\textcolor{black}{Next, to state a corollary of this proposition, we introduce some more notation.}
For $F_1,F_2,F_3,F_4\in L_\eta^2$ such that $\EE[F_i]=0$, $i=1,2,3,4$, we set
\begin{align}
\Gamma_1(F_1)^2:&=
4\int_{\mathbb{X}^3}\mathbb{E}[(D_{x_1,x_3}^{2}F_1)^4]^{1/4}\mathbb{E}[(D_{x_2,x_3}^{2}F_1)^4]^{1/4}\mathbb{E}[(D_{x_1}F_1)^4]^{1/4}\mathbb{E}[(D_{x_2}F_1)^4]^{1/4}\mu^{\otimes 3}(\mathrm{d}x_1,\dd x_2,\dd x_3)\nonumber\\
&\qquad\qquad\qquad
+\int_{\mathbb{X}^3}\mathbb{E}[(D_{x_1,x_3}^{2}F_1)^4]^{1/2}\mathbb{E}[(D_{x_2,x_3}^{2}F_1)^4]^{1/2}\mu^{\otimes 3}(\mathrm{d}x_1,\dd x_2,\dd x_3)\nonumber\\
&=
4\int_{\mathbb{X}}\left(\int_{\mathbb{X}}\mathbb{E}[(D_{x_1,x_2}^{2}F_1)^4]^{1/4}\mathbb{E}[(D_{x_1}F_1)^4]^{1/4}\mu(\dd x_1)\right)^2\mu(\dd x_2)\nonumber\\
&\qquad\qquad\qquad\qquad\qquad\qquad
+\int_{\mathbb{X}}\left(\int_{\mathbb X}\mathbb{E}[(D_{x_1,x_2}^{2}F_1)^4]^{1/2}\mu(\dd x_1)\right)^2\mu(\dd x_2),\nonumber
\end{align}
\begin{align}
\Gamma(F_1,F_2,F_3,F_4)^2:&=4\Lambda(F_1,F_2,F_3,F_4)^2
+\int_{\mathbb{X}^3}\mathbb{E}[(D_{x_1,x_3}^{2}F_1)^4]^{1/2}\mathbb{E}[(D_{x_2,x_3}^{2}F_2)^4]^{1/2}\mu^{\otimes 3}(\mathrm{d}x_1,\dd x_2,\dd x_3)\nonumber\\
&=4\Lambda(F_1,F_2,F_3,F_4)^2\nonumber\\
&\qquad\qquad
+\int_{\mathbb{X}}\left(\int_{\mathbb X}\mathbb{E}[(D_{x_1,x_2}^{2}F_1)^4]^{1/2}\mu(\dd x_1)\right)
\left(\int_{\mathbb X}\mathbb{E}[(D_{x_1,x_2}^{2}F_2)^4]^{1/2}\mu(\dd x_1)\right)\mu(\dd x_2),\nonumber
\end{align}
\begin{align}
&\Lambda(F_1,F_2,F_3,F_4)^2\nonumber\\
:&=\int_{\mathbb{X}^3}\mathbb{E}[(D_{x_1,x_3}^{2}F_1)^4]^{1/4}\mathbb{E}[(D_{x_2,x_3}^{2}F_2)^4]^{1/4}\mathbb{E}[(D_{x_1}F_3)^4]^{1/4}\mathbb{E}[(D_{x_2}F_4)^4]^{1/4}\mu^{\otimes 3}(\mathrm{d}x_1,\dd x_2,\dd x_3)\nonumber\\
&=\int_{\mathbb{X}}\left(\int_{\mathbb{X}}\mathbb{E}[(D_{x_1,x_2}^{2}F_1)^4]^{1/4}
\mathbb{E}[(D_{x_1}F_3)^4]^{1/4}\mu(\dd x_1)\right)
\left(\int_{\mathbb{X}}\mathbb{E}[(D_{x_1,x_2}^{2}F_2)^4]^{1/4}
\mathbb{E}[(D_{x_1}F_4)^4]^{1/4}\mu(\dd x_1)\right)\mu(\dd x_2)\nonumber
\end{align}
and
\[
\Theta(F_1,F_2):=\int_{\mathbb X}\mathbb{E}[(D_x F_1)^4]^{1/4}\EE[(D_x F_2)^4]^{1/4}\mu(\mathrm{d}x).
\]
Note that $\Gamma_1(F_1)^2=\Gamma(F_1,F_1,F_1,F_1)^2$,
\begin{equation}\label{eq:fubinibis}
\Lambda(F_1,F_2,F_3,F_4)^2=\Lambda(F_2,F_1,F_4,F_3)^2\quad\text{and}\quad\Gamma(F_1,F_2,F_3,F_4)^2=\Gamma(F_2,F_1,F_4,F_3)^2.
\end{equation}
Note also that by the Cauchy-Schwarz inequality
\begin{equation}\label{eq:gammalambdaless}
\gamma(F_1,F_2,F_3,F_4)^2\leq\Gamma(F_1,F_2,F_3,F_4)^2,\quad\lambda(F_1,F_2,F_3,F_4)^2\leq\Lambda(F_1,F_2,F_3,F_4)^2,
\quad\theta(F_1,F_2)\leq\Theta(F_1,F_2)
\end{equation}
(and, so, in particular, $\gamma_1(F_1)\leq\Gamma_1(F_1)$).

The following corollary is an immediate consequence of the inequalities \eqref{eq:gammalambdaless} and Proposition \ref{thm:ASCLT}.

\begin{Corollary}\label{cor:ASCLT}
Assume $\{F_n\}_{n\geq 1}\subset L_\eta^2$ and $\mathbb{E}[F_n]=0$, $\mathbb{E}[F_n^2]=1$, for any $n\geq 1$. Moreover, suppose
\eqref{hyp:convgamma} and \eqref{hyp:series1conv} with $\Gamma_1$ in place of $\gamma_1$,
\eqref{hyp:series2conv} with $\Theta$ in place of $\theta$, \eqref{hyp:series3conv} with $\Gamma$ in place of $\gamma$,
and \eqref{hyp:series4conv} with $\Lambda$ in place of $\lambda$. Then $\{F_n\}_{n\geq 1}$ satisfies the ASCLT.
\end{Corollary}

\noindent{\it Proof\,\,of\,\,Theorem\,\,\ref{cor:ASCLTAPP}.} For $n\geq 1$, $x_1,x_2\in\mathbb X$, we have
\[
D_{x_1}F_n=\frac{D_{x_1}H_n}{\sqrt{\mathbb{V}ar(H_n)}}\quad\text{and}\quad
D_{x_1,x_2}^2 F_n=\frac{D_{x_1,x_2}^2 H_n}{\sqrt{\mathbb{V}ar(H_n)}}.
\]
So, letting $C\in\mathbb R_+$ denote a positive constant which may vary from line to line,
by H\"older's inequality and assumptions \eqref{eq:os1app} and \eqref{eq:os2app}, for any $n\geq 1$, we have
\begin{align}
\mathbb{E}[|D_{x_1}F_n|^4]
%&=\mathbb{E}[\bold{1}\{D_{x_1}F_n\neq 0\}|D_{x_1}F_n|^4]
%\nonumber\\
%&\leq\mathbb{E}[\bold{1}\{D_{x_1}F_n\neq 0\}^{(4+p)/p}]^{p/(4+p)}\mathbb{E}[|D_{x_1}F_n|^{4+p}]^{4/(4+p)}\nonumber\\
%&=\mathbb{P}(D_{x_1}F_n\neq 0)^{p/(4+p)}\mathbb{E}[|D_{x_1}F_n|^{4+p}]^{4/(4+p)}\nonumber\\
&\leq\mathbb{P}(D_{x_1}H_n\neq 0)^{p/(4+p)}\frac{\mathbb{E}[|D_{x_1}H_n|^{4+p}]^{4/(4+p)}}{\mathbb{V}ar(H_n)^2}\nonumber\\
&\leq C\frac{\mathbb{P}(D_{x_1}H_n\neq 0)^{p/(4+p)}}{\mathbb{V}ar(H_n)^2},\quad\text{$\mu$-a.e. $x_1\in\mathbb X$}\label{eq:os4}
\end{align}
\begin{align}
\mathbb{E}[|D_{x_1,x_2}^2 F_n|^4]
%&\leq
%\mathbb{P}(D_{x_1,x_2}^2 F_n\neq 0)^{q/(4+q)}\mathbb{E}[|D_{x_1,x_2}^2 F_n|^{4+q}]^{4/(4+q)}\nonumber\\
%&=\mathbb{P}(D_{x_1,x_2}^2 H_n\neq 0)^{q/(4+q)}\frac{\mathbb{E}[|D_{x_1,x_2}H_n|^{4+q}]^{4/(4+q)}}{\mathbb{V}ar(H_n)^2}\nonumber\\
&\leq C\frac{\mathbb{P}(D_{x_1,x_2}^2 H_n\neq 0)^{q/(4+q)}}{\mathbb{V}ar(H_n)^2},\quad\text{$\mu^{\otimes 2}$-a.e. $(x_1,x_2)\in\mathbb X^2$}\label{eq:os5}
\end{align}
and
\begin{align}
\mathbb{E}[|D_{x_1}F_n|^3]
%&=\mathbb{E}[\bold{1}\{D_{x_1}F_n\neq 0\}|D_{x_1}F_n|^3]
%\nonumber\\
%&\leq\mathbb{E}[\bold{1}\{D_{x_1}F_n\neq 0\}^{(4+p)/(1+p)}]^{(1+p)/(4+p)}\mathbb{E}[|D_{x_1}F_n|^{4+p}]^{3/(4+p)}\nonumber\\
%&=\mathbb{P}(D_{x_1}F_n\neq 0)^{(1+p)/(4+p)}\mathbb{E}[|D_{x_1}F_n|^{4+p}]^{3/(4+p)}\nonumber\\
&\leq\mathbb{P}(D_{x_1}H_n\neq 0)^{(1+p)/(4+p)}\frac{\mathbb{E}[|D_{x_1}H_n|^{4+p}]^{3/(4+p)}}{\mathbb{V}ar(H_n)^{3/2}}\nonumber\\
&\leq C\frac{\mathbb{P}(D_{x_1}H_n\neq 0)^{(1+p)/(4+p)}}{\mathbb{V}ar(H_n)^{3/2}},\quad\text{$\mu$-a.e. $x_1\in\mathbb X$.}\label{eq:osmomter}
\end{align}
For $k,l,m,n\geq 1$, the inequalities \eqref{eq:os4}, \eqref{eq:os5} and \eqref{eq:osmomter} imply
\[
\Gamma_1(F_k)^2\leq C\bold{\Gamma}_1(H_k)^2,\quad\Gamma(F_k,F_l,F_m,F_n)^2\leq C\bold{\Gamma}(H_k,H_l,H_m,H_n)^2,
\]
\[
\Lambda(F_k,F_l,F_m,F_n)^2\leq C\bold{\Lambda}(H_k,H_l,H_m,H_n)^2,\quad
\Theta(F_k,F_l)^2\leq C\bold{\Theta}(H_k,H_l)^2,\quad\gamma_2(F_k)\leq C\bold{\Gamma}_2(H_k).
\]
The claim follows by these relations and Corollary \ref{cor:ASCLT}.
\\
\noindent$\square$

\noindent{\it Proof\,\,of\,\,Proposition\,\,\ref{thm:ASCLT}.} We start noticing that in fact $F_n\in\mathrm{Dom}(D)$, $n\geq 1$. Indeed,
by \eqref{hyp:series2conv} it follows
\[
\int_{\mathbb X}\mathbb{E}[|D_x F|^2]\mu(\dd x)=\theta(F_n,F_n)<\infty,\quad n\geq 1.
\]
To prove the ASCLT we are going to apply the Ibragimov and Lifshits criterion stated as Lemma \ref{thm:IL}.
So, let $\Delta_n(t)$ be given by \eqref{eq:deltan} with $F_k$ in place of $X_k$ and $Z$ in place of $X$. We have
\begin{align}
\mathbb{E}[|\Delta_n(t)|^2]&=\frac{1}{(\log n)^2}\sum_{k,l}^{1,n}\frac{1}{kl}\EE\left[\left(\mathrm{e}^{\bold{i}tF_k}-\mathrm{e}^{-t^2/2}\right)\left(\mathrm{e}^{-\bold{i}tF_l}-\mathrm{e}^{-t^2/2}\right)\right]\nonumber\\
&=\frac{1}{(\log n)^2}\sum_{k,l}^{1,n}\frac{1}{kl}\Biggl[\left(\EE[\mathrm{e}^{\bold{i}t(F_k-F_l)}]-\mathrm{e}^{-t^2}\right)
-\mathrm{e}^{-t^2/2}(\EE[\mathrm{e}^{\bold{i}tF_k}]-\mathrm{e}^{-t^2/2})\nonumber\\
&\qquad\qquad\qquad\qquad\qquad\qquad
-\mathrm{e}^{-t^2/2}\left(\EE[\mathrm{e}^{-\bold{i}t F_l}]-\mathrm{e}^{-t^2/2}\right)\Biggr].\nonumber
%\\
%&=\frac{1}{(\log n)^2}\sum_{k,l}^{1,n}\frac{1}{kl}\Biggl[\left(\varphi_{F_k-F_l}(t)-\mathrm{e}^{-t^2}\right)-\mathrm{e}^{-t^2/2}(\varphi_{F_k}(t)-\mathrm{e}^{-t^2/2})
%-\mathrm{e}^{-t^2/2}\left(\varphi_{F_l}(-t)-\mathrm{e}^{-t^2/2}\right)\Biggr].\nonumber
\end{align}
Therefore, for $r\in\mathbb R_+$,
\begin{align}
\sup_{|t|\leq r}\sum_{n\geq 2}\frac{\mathbb{E}[|\Delta_n(t)|^2]}{n\log n}&\leq\sup_{|t|\leq r}\sum_{n\geq 2}\frac{1}{n(\log n)^3}\sum_{k,l}^{1,n}\frac{|\EE[\mathrm{e}^{\bold{i}t(F_k-F_l)}]-\mathrm{e}^{-t^2}|}{kl}
\label{eq:1stsup}\\
&\quad
+\sup_{|t|\leq r}\sum_{n\geq 2}\frac{1}{n(\log n)^3}\sum_{l=1}^{n}\frac{1}{l}\sum_{k=1}^{n}\frac{|\EE[\mathrm{e}^{\bold{i}tF_k}]-\mathrm{e}^{-t^2/2}|}{k}\label{eq:2ndsup}\\
&\quad
+\sup_{|t|\leq r}\sum_{n\geq 2}\frac{1}{n(\log n)^3}\sum_{k=1}^{n}\frac{1}{k}\sum_{l=1}^{n}\frac{|\EE[\mathrm{e}^{-\bold{i}tF_l}]-\mathrm{e}^{-t^2/2}|}{l}.\label{eq:3rdsup}
\end{align}
Let $t\in\mathbb R$ and $r\in\mathbb{R}_+$ be such that $|t|\leq r$. \textcolor{black}{Since $\mathbb{E}[F_n^2]=1$ for each $n$,} by Lemma \ref{le:ledisfcar} we have
\[
|\EE[\mathrm{e}^{\bold{i}F_k t}]-\mathrm{e}^{-t^2/2}|\leq r^2\gamma_1(F_k)+\frac{r^3}{\sqrt 2}\gamma_2(F_k).
\]
Therefore,
\begin{align}
\sup_{|t|\leq r}\sum_{n\geq 2}\frac{1}{n(\log n)^3}\sum_{l=1}^{n}\frac{1}{l}\sum_{k=1}^{n}\frac{|\EE[\mathrm{e}^{\bold{i}tF_k}]-\mathrm{e}^{-t^2/2}|}{k}
\leq\sum_{n\geq 2}\frac{1}{n(\log n)^3}\sum_{l=1}^{n}\frac{1}{l}\sum_{k=1}^{n}\frac{r^2\gamma_1(F_k)+2^{-1/2}r^3\gamma_2(F_k)}{k}.\nonumber
\end{align}
Since $\lim_{n\to\infty}\sum_{l=1}^{n}l^{-1}/\log n=1$, the finiteness of the term in \eqref{eq:2ndsup} is guaranteed by the assumption \eqref{hyp:series1conv}.
Since by Lemma \ref{le:ledisfcar} we also have
\[
|\EE[\mathrm{e}^{-\bold{i}tF_l}]-\mathrm{e}^{-t^2/2}|\leq r^2\gamma_1(F_l)+\frac{r^3}{\sqrt 2}\gamma_2(F_l),
\]
the same argument guarantees the finiteness of the term in \eqref{eq:3rdsup}. As far as the finiteness of the term
in the right-hand side of \eqref{eq:1stsup}, we start noticing that, again by Lemma \ref{le:ledisfcar} we get
\begin{align}
|\EE[\mathrm{e}^{\bold{i}t(F_k-F_l)}]-\mathrm{e}^{-t^2}|&=|\EE[\mathrm{e}^{\bold{i}(\sqrt{2}t)(F_k-F_l)/\sqrt 2}]-\mathrm{e}^{-(\sqrt{2}t)^2/2}|\nonumber\\
&\leq 2|t|^2|\mathbb{E}[F_kF_l]|+2|t|^2\gamma_1\left(\frac{F_k-F_l}{\sqrt 2}\right)+2|t|^3\gamma_2\left(\frac{F_k-F_l}{\sqrt 2}\right)\nonumber\\
&\leq 2\textcolor{black}{r}^2\mathbb{E}[F_kF_l]|+2\textcolor{black}{r}^2\gamma_1\left(\frac{F_k-F_l}{\sqrt 2}\right)+2\textcolor{black}{r}^3\gamma_2\left(\frac{F_k-F_l}{\sqrt 2}\right).\nonumber
\end{align}
Therefore the term in the right-hand side of \eqref{eq:1stsup} is finite if
\begin{align}
&\sum_{n\geq 2}\frac{1}{n(\log n)^3}\sum_{k,l}^{1,n}\frac{|\EE[F_k F_l]|}{kl}<\infty,\label{eq:series1}\\
&\sum_{n\geq 2}\frac{1}{n(\log n)^3}\sum_{k,l}^{1,n}\frac{\gamma_2\left(\frac{F_k-F_l}{\sqrt 2}\right)}{kl}<\infty,\label{eq:series2}\\
&\sum_{n\geq 2}\frac{1}{n(\log n)^3}\sum_{k,l}^{1,n}\frac{\gamma_1\left(\frac{F_k-F_l}{\sqrt 2}\right)}{kl}<\infty.\label{eq:series3}
\end{align}
\noindent\textcolor{black}{\it{Proof\,\,of\,\,\eqref{eq:series1}.}}
Due to the symmetry (with respect to $k$ and $l$) of the mapping $(k,l)\mapsto\frac{|\EE[F_k F_l]|}{kl}$, it holds
\begin{align}
\sum_{n\geq 2}\frac{1}{n(\log n)^3}\sum_{k,l}^{1,n}\frac{|\EE[F_k F_l]|}{kl}&
=\sum_{n\geq 2}\frac{1}{n(\log n)^3}\sum_{k=1}^{n}\frac{1}{k^2}+2\sum_{n\geq 2}\frac{1}{n(\log n)^3}\sum_{l=2}^{n}\sum_{k=1}^{l-1}\frac{|\EE[F_k F_l]|}{kl}.\nonumber
\end{align}
Since $\sum_{n\geq 2}\frac{1}{n(\log n)^\beta}<\infty$, $\beta>1$, and $\sum_{n\geq 2}\frac{1}{n^2}<\infty$, the infinite sum
\eqref{eq:series1} converges if and only if
\begin{equation}\label{eq:series1bis}
\sum_{n\geq 2}\frac{1}{n(\log n)^3}\sum_{l=1}^{n}\sum_{k=1}^{l}\frac{|\EE[F_k F_l]|}{kl}<\infty.
\end{equation}
By Lemmas \ref{le:covariance}, \ref{le:DLinverse} and the Cauchy-Schwarz inequality, we get
\begin{align}
|\EE[F_k F_l]|&\leq\EE[|\langle DF_k,-DL^{-1}F_l\rangle_1|]\nonumber\\
&\leq\int_{\mathbb X}\EE[|D_x F_k||D_x L^{-1}F_l|]\mu(\dd x)\leq\theta(F_k,F_l),\nonumber
\end{align}
and therefore \eqref{eq:series1bis} is a consequence of the assumption \eqref{hyp:series2conv}.
\\
\noindent\textcolor{black}{\it{Proof\,\,of\,\,\eqref{eq:series2}.}}
By the inequality $(a+b)^3\leq 4(a^3+b^3)$, $a,b\geq 0$, we have
\begin{align}
\gamma_2\left(\frac{F_k-F_l}{\sqrt 2}\right)&=2^{-3/2}\int_{\mathbb{X}}\mathbb{E}[|D_xF_k-D_xF_l|^3]\mu(\mathrm{d}x)\nonumber\\
%&\leq\int_{\mathbb{X}}\mathbb{E}[|D_xF_k|^3]\mu(\mathrm{d}x)+\int_{\mathbb{X}}\mathbb{E}[|D_xF_l|^3]\mu(\mathrm{d}x)\nonumber\\
&\leq\sqrt{2}(\gamma_2(F_k)+\gamma_2(F_l))\nonumber.
\end{align}
Therefore, \eqref{eq:series2} follows by the assumption \eqref{hyp:series1conv} (with $i=2$).\\
\textcolor{black}{\it{Proof\,\,of\,\,\eqref{eq:series3}.}} Due to the symmetry (with respect to $k$ and $l$) of the mapping
$(k,l)\mapsto\frac{\gamma_{1}\left(\frac{F_k-F_l}{\sqrt 2}\right)}{kl}$, the infinite sum \eqref{eq:series3} converges if and only if
\begin{equation}\label{eq:series2bis}
\sum_{n\geq 2}\frac{1}{n(\log n)^3}\sum_{l=1}^{n}\sum_{k=1}^{l}\frac{\gamma_1\left(\frac{F_k-F_l}{\sqrt 2}\right)}{kl}<\infty.
\end{equation}
Using the inequalities $(a+b)^2\leq 2a^2+2b^2$ and $(a+b)^{1/2}\leq a^{1/2}+b^{1/2}$, $a,b\geq 0$,
we now bound the quantity $\gamma_1((F_k-F_l)/\sqrt 2)$. Hereon, for simplicity of notation we put
$\mu^{\otimes 3}(\mathrm{d}\bold x):=\mu^{\otimes 3}(\dd x_1,\dd x_2, \dd x_3)$. We have:
\begin{align}
&\gamma_1\left(\frac{F_k-F_l}{\sqrt 2}\right)\leq
\Biggl(
4\int_{\mathbb{X}^3}(\mathbb{E}[(|D_{x_1,x_3}^{2}F_k|^2+|D_{x_1,x_3}^{2}F_l|^2)(|D_{x_2,x_3}^{2}F_k|^2+|D_{x_2,x_3}^{2}F_l|^2)])^{1/2}\nonumber\\
&\qquad\qquad\qquad\qquad
\times(\mathbb{E}[(|D_{x_1}F_k|^2+|D_{x_1}F_l|^2)(|D_{x_2}F_k|^2+|D_{x_2}F_l|^2)])^{1/2}\mu^{\otimes 3}(\mathrm{d}\bold{x})\nonumber\\
&\qquad
+\int_{\mathbb{X}^3}\mathbb{E}[(|D_{x_1,x_3}^{2}F_k|^2+|D_{x_1,x_3}^{2}F_l|^2)(|D_{x_2,x_3}^{2}F_k|^2+|D_{x_2,x_3}^{2}F_l|^2)]
\mu^{\otimes 3}(\mathrm{d}\bold{x})\Biggr)^{1/2}\nonumber\\
&=\Biggl(
4\int_{\mathbb{X}^3}\Biggl(
(\mathbb{E}[|D_{x_1,x_3}^{2}F_k|^2|D_{x_2,x_3}^{2}F_k|^2]+\mathbb{E}[|D_{x_1,x_3}^{2}F_k|^2|D_{x_2,x_3}^{2}F_l|^2]\nonumber\\
&\qquad\qquad\qquad\qquad
+\mathbb{E}[|D_{x_1,x_3}^{2}F_l|^2|D_{x_2,x_3}^{2}F_k|^2]+\mathbb{E}[|D_{x_1,x_3}^{2}F_l|^2|D_{x_2,x_3}^{2}F_l|^2])\nonumber\\
&\qquad\qquad
\times
(\mathbb{E}[|D_{x_1}F_k|^2|D_{x_2}F_k|^2]+\mathbb{E}[|D_{x_1}F_k|^2|D_{x_2}F_l|^2]\nonumber\\
&\qquad\qquad\qquad\qquad
+\mathbb{E}[|D_{x_1}F_l|^2|D_{x_2}F_k|^2]+\mathbb{E}[|D_{x_1}F_l|^2|D_{x_2}F_l|^2])\Biggr)^{1/2}\mu^{\otimes 3}(\mathrm{d}\bold x)\nonumber\\
&\qquad
+\int_{\mathbb{X}^3}(\mathbb{E}[|D_{x_1,x_3}^{2}F_k|^2|D_{x_2,x_3}^{2}F_k|^2]+\mathbb{E}[|D_{x_1,x_3}^{2}F_k|^2|D_{x_2,x_3}^{2}F_l|^2]\nonumber\\
&\qquad\qquad\qquad\qquad
+\mathbb{E}[|D_{x_1,x_3}^{2}F_l|^2|D_{x_2,x_3}^{2}F_k|^2]+\mathbb{E}[
|D_{x_1,x_3}^{2}F_l|^2|D_{x_2,x_3}^{2}F_l|^2])\mu^{\otimes 3}(\mathrm{d}\bold x)\Biggr)^{1/2}\nonumber\\
&=\Biggl(
4\int_{\mathbb{X}^3}\mu^{\otimes 3}(\mathrm{d}\bold x)\nonumber\\
&\Biggl(\mathbb{E}[|D_{x_1,x_3}^{2}F_k|^2|D_{x_2,x_3}^{2}F_k|^2]\mathbb{E}[|D_{x_1}F_k|^2|D_{x_2}F_k|^2]
+\mathbb{E}[|D_{x_1,x_3}^{2}F_k|^2|D_{x_2,x_3}^{2}F_k|^2]\mathbb{E}[|D_{x_1}F_k|^2|D_{x_2}F_l|^2]\nonumber\\
&\qquad
+\mathbb{E}[|D_{x_1,x_3}^{2}F_k|^2|D_{x_2,x_3}^{2}F_k|^2]\mathbb{E}[|D_{x_1}F_l|^2|D_{x_2}F_k|^2]
+\mathbb{E}[|D_{x_1,x_3}^{2}F_k|^2|D_{x_2,x_3}^{2}F_k|^2]\mathbb{E}[|D_{x_1}F_l|^2|D_{x_2}F_l|^2]\nonumber\\
&\qquad
+\mathbb{E}[|D_{x_1,x_3}^{2}F_k|^2|D_{x_2,x_3}^{2}F_l|^2]\mathbb{E}[|D_{x_1}F_k|^2|D_{x_2}F_k|^2]
+\mathbb{E}[|D_{x_1,x_3}^{2}F_k|^2|D_{x_2,x_3}^{2}F_l|^2]\mathbb{E}[|D_{x_1}F_k|^2|D_{x_2}F_l|^2]\nonumber\\
&\qquad
+\mathbb{E}[|D_{x_1,x_3}^{2}F_k|^2|D_{x_2,x_3}^{2}F_l|^2]\mathbb{E}[|D_{x_1}F_l|^2|D_{x_2}F_k|^2]
+\mathbb{E}[|D_{x_1,x_3}^{2}F_k|^2|D_{x_2,x_3}^{2}F_l|^2]\mathbb{E}[|D_{x_1}F_l|^2|D_{x_2}F_l|^2]\nonumber\\
&\qquad
+\mathbb{E}[|D_{x_1,x_3}^{2}F_l|^2|D_{x_2,x_3}^{2}F_k|^2]\mathbb{E}[|D_{x_1}F_k|^2|D_{x_2}F_k|^2]
+\mathbb{E}[|D_{x_1,x_3}^{2}F_l|^2|D_{x_2,x_3}^{2}F_k|^2]\mathbb{E}[|D_{x_1}F_k|^2|D_{x_2}F_l|^2]\nonumber\\
&\qquad
+\mathbb{E}[|D_{x_1,x_3}^{2}F_l|^2|D_{x_2,x_3}^{2}F_k|^2]\mathbb{E}[|D_{x_1}F_l|^2|D_{x_2}F_k|^2]
+\mathbb{E}[|D_{x_1,x_3}^{2}F_l|^2|D_{x_2,x_3}^{2}F_k|^2]\mathbb{E}[|D_{x_1}F_l|^2|D_{x_2}F_l|^2]\nonumber\\
&\qquad
+\mathbb{E}[|D_{x_1,x_3}^{2}F_l|^2|D_{x_2,x_3}^{2}F_l|^2]\mathbb{E}[|D_{x_1}F_k|^2|D_{x_2}F_k|^2]
+\mathbb{E}[|D_{x_1,x_3}^{2}F_l|^2|D_{x_2,x_3}^{2}F_l|^2]\mathbb{E}[|D_{x_1}F_k|^2|D_{x_2}F_l|^2]\nonumber\\
&\qquad
+\mathbb{E}[|D_{x_1,x_3}^{2}F_l|^2|D_{x_2,x_3}^{2}F_l|^2]\mathbb{E}[|D_{x_1}F_l|^2|D_{x_2}F_k|^2]
+\mathbb{E}[|D_{x_1,x_3}^{2}F_l|^2|D_{x_2,x_3}^{2}F_l|^2]\mathbb{E}[|D_{x_1}F_l|^2|D_{x_2}F_l|^2]\Biggr)^{1/2}\nonumber\\
&\qquad
+\int_{\mathbb{X}^3}(\mathbb{E}[|D_{x_1,x_3}^{2}F_k|^2|D_{x_2,x_3}^{2}F_k|^2]+\mathbb{E}[|D_{x_1,x_3}^{2}F_k|^2|D_{x_2,x_3}^{2}F_l|^2]\nonumber\\
&\qquad\qquad\qquad\qquad
+\mathbb{E}[|D_{x_1,x_3}^{2}F_l|^2|D_{x_2,x_3}^{2}F_k|^2]+\mathbb{E}[
|D_{x_1,x_3}^{2}F_l|^2|D_{x_2,x_3}^{2}F_l|^2])\mu^{\otimes 3}(\mathrm{d}\bold x)\Biggr)^{1/2}\nonumber\\
&\leq\Biggl[\int_{\mathbb{X}^3}(4(\mathbb{E}[|D_{x_1,x_3}^{2}F_k|^2|D_{x_2,x_3}^{2}F_k|^2])^{1/2}(\mathbb{E}[|D_{x_1}F_k|^2|D_{x_2}F_k|^2])^{1/2}
+\mathbb{E}[|D_{x_1,x_3}^{2}F_k|^2|D_{x_2,x_3}^{2}F_k|^2])\mu^{\otimes 3}(\mathrm{d}\bold x)\nonumber\\
&+\int_{\mathbb{X}^3}(4(\mathbb{E}[|D_{x_1,x_3}^{2}F_l|^2|D_{x_2,x_3}^{2}F_l|^2])^{1/2}(\mathbb{E}[|D_{x_1}F_l|^2|D_{x_2}F_l|^2])^{1/2}
+\mathbb{E}[|D_{x_1,x_3}^{2}F_l|^2|D_{x_2,x_3}^{2}F_l|^2])\mu^{\otimes 3}(\mathrm{d}\bold x)\nonumber\\
&+\int_{\mathbb{X}^3}(4(\mathbb{E}[|D_{x_1,x_3}^{2}F_k|^2|D_{x_2,x_3}^{2}F_l|^2])^{1/2}(\mathbb{E}[|D_{x_1}F_k|^2|D_{x_2}F_l|^2])^{1/2}
+\mathbb{E}[|D_{x_1,x_3}^{2}F_k|^2|D_{x_2,x_3}^{2}F_l|^2])\mu^{\otimes 3}(\mathrm{d}\bold x)\nonumber\\
&+\int_{\mathbb{X}^3}(4(\mathbb{E}[|D_{x_1,x_3}^{2}F_l|^2|D_{x_2,x_3}^{2}F_k|^2])^{1/2}(\mathbb{E}[|D_{x_1}F_l|^2|D_{x_2}F_k|^2])^{1/2}
+\mathbb{E}[|D_{x_1,x_3}^{2}F_l|^2|D_{x_2,x_3}^{2}F_k|^2])\mu^{\otimes 3}(\mathrm{d}\bold x)\nonumber\\
&+4\int_{\mathbb X^3}(\mathbb{E}[|D_{x_1,x_3}^{2}F_k|^2|D_{x_2,x_3}^{2}F_k|^2])^{1/2}(\mathbb{E}[|D_{x_1}F_k|^2|D_{x_2}F_l|^2])^{1/2}\mu^{\otimes 3}(\mathrm{d}\mathbf x)\nonumber\\
&+4\int_{\mathbb X^3}(\mathbb{E}[|D_{x_1,x_3}^{2}F_k|^2|D_{x_2,x_3}^{2}F_k|^2])^{1/2}(\mathbb{E}[|D_{x_1}F_l|^2|D_{x_2}F_k|^2])^{1/2}\mu^{\otimes 3}(\mathrm{d}\mathbf x)\nonumber\\
&+4\int_{\mathbb X^3}(\mathbb{E}[|D_{x_1,x_3}^{2}F_k|^2|D_{x_2,x_3}^{2}F_k|^2])^{1/2}(\mathbb{E}[|D_{x_1}F_l|^2|D_{x_2}F_l|^2])^{1/2}\mu^{\otimes 3}(\mathrm{d}\mathbf x)\nonumber\\
&+4\int_{\mathbb X^3}(\mathbb{E}[|D_{x_1,x_3}^{2}F_k|^2|D_{x_2,x_3}^{2}F_l|^2])^{1/2}(\mathbb{E}[|D_{x_1}F_k|^2|D_{x_2}F_k|^2])^{1/2}\mu^{\otimes 3}(\mathrm{d}\mathbf x)\nonumber\\
&+4\int_{\mathbb X^3}(\mathbb{E}[|D_{x_1,x_3}^{2}F_k|^2|D_{x_2,x_3}^{2}F_l|^2])^{1/2}(\mathbb{E}[|D_{x_1}F_l|^2|D_{x_2}F_k|^2])^{1/2}\mu^{\otimes 3}(\mathrm{d}\mathbf x)\nonumber\\
&+4\int_{\mathbb X^3}(\mathbb{E}[|D_{x_1,x_3}^{2}F_k|^2|D_{x_2,x_3}^{2}F_l|^2])^{1/2}(\mathbb{E}[|D_{x_1}F_l|^2|D_{x_2}F_l|^2])^{1/2}\mu^{\otimes 3}(\mathrm{d}\mathbf x)\nonumber\\
&+4\int_{\mathbb X^3}(\mathbb{E}[|D_{x_1,x_3}^{2}F_l|^2|D_{x_2,x_3}^{2}F_k|^2])^{1/2}(\mathbb{E}[|D_{x_1}F_k|^2|D_{x_2}F_k|^2])^{1/2}\mu^{\otimes 3}(\mathrm{d}\mathbf x)\nonumber\\
&+4\int_{\mathbb X^3}(\mathbb{E}[|D_{x_1,x_3}^{2}F_l|^2|D_{x_2,x_3}^{2}F_k|^2])^{1/2}(\mathbb{E}[|D_{x_1}F_k|^2|D_{x_2}F_l|^2])^{1/2}\mu^{\otimes 3}(\mathrm{d}\mathbf x)\nonumber\\
&+4\int_{\mathbb X^3}(\mathbb{E}[|D_{x_1,x_3}^{2}F_l|^2|D_{x_2,x_3}^{2}F_k|^2])^{1/2}(\mathbb{E}[|D_{x_1}F_l|^2|D_{x_2}F_l|^2])^{1/2}\mu^{\otimes 3}(\mathrm{d}\mathbf x)\nonumber\\
&+4\int_{\mathbb X^3}(\mathbb{E}[|D_{x_1,x_3}^{2}F_l|^2|D_{x_2,x_3}^{2}F_l|^2])^{1/2}(\mathbb{E}[|D_{x_1}F_k|^2|D_{x_2}F_k|^2])^{1/2}\mu^{\otimes 3}(\mathrm{d}\mathbf x)\nonumber\\
&+4\int_{\mathbb X^3}(\mathbb{E}[|D_{x_1,x_3}^{2}F_l|^2|D_{x_2,x_3}^{2}F_l|^2])^{1/2}(\mathbb{E}[|D_{x_1}F_k|^2|D_{x_2}F_l|^2])^{1/2}\mu^{\otimes 3}(\mathrm{d}\mathbf x)\nonumber\\
&+4\int_{\mathbb X^3}(\mathbb{E}[|D_{x_1,x_3}^{2}F_l|^2|D_{x_2,x_3}^{2}F_l|^2])^{1/2}(\mathbb{E}[|D_{x_1}F_l|^2|D_{x_2}F_k|^2])^{1/2}\mu^{\otimes 3}(\mathrm{d}\mathbf x)\Biggr]^{1/2}.\label{eq:ineqbig}
\end{align}
Finally, by \eqref{eq:ineqbig}, using firstly again the inequality $(a+b)^{1/2}\leq a^{1/2}+b^{1/2}$, $a,b\geq 0$, and secondly the identities \eqref{eq:fubini},
we have
\begin{align}
\gamma_1\left(\frac{F_k-F_l}{\sqrt 2}\right)&\leq\gamma_1(F_k)+\gamma_1(F_l)+\gamma(F_k,F_l,F_k,F_l)+\gamma(F_l,F_k,F_l,F_k)\nonumber\\
&+\lambda(F_k,F_k,F_k,F_l)+\lambda(F_k,F_k,F_l,F_k)+\lambda(F_k,F_k,F_l,F_l)+\lambda(F_k,F_l,F_k,F_k)\nonumber\\
&+\lambda(F_k,F_l,F_l,F_k)+\lambda(F_k,F_l,F_l,F_l)+\lambda(F_l,F_k,F_k,F_k)+\lambda(F_l,F_k,F_k,F_l)\nonumber\\
&+\lambda(F_l,F_k,F_l,F_l)+\lambda(F_l,F_l,F_k,F_k)+\lambda(F_l,F_l,F_k,F_l)+\lambda(F_l,F_l,F_l,F_k)\nonumber\\
&=\gamma_1(F_k)+\gamma_1(F_l)+2\gamma(F_k,F_l,F_k,F_l)+2\lambda(F_k,F_k,F_k,F_l)\nonumber\\
&+\lambda(F_k,F_k,F_l,F_l)+2\lambda(F_l,F_k,F_k,F_k)+2\lambda(F_k,F_l,F_l,F_k)+2\lambda(F_k,F_l,F_l,F_l)\nonumber\\
&+\lambda(F_l,F_l,F_k,F_k)+2\lambda(F_l,F_l,F_l,F_k).\nonumber
\end{align}
Therefore \eqref{eq:series2bis} follows by the assumptions \eqref{hyp:series1conv} (with $i=1$), \eqref{hyp:series3conv} and \eqref{hyp:series4conv}.

By Lemma \ref{le:ledisfcar}, assumption \eqref{hyp:convgamma} and L\'evy's continuity theorem, the sequence $\{F_n\}_{n\geq 1}$ converges in law to $Z$. At last, the claim follows by the Ibragimov and Lifshits criterion.
\\
\noindent$\square$
\\
\\
$\bold{Acknowledgments}$ We would like to thank Matthias Schulte for useful discussions about his paper \cite{LSY}, and
\textcolor{black}{the anonymous referees for a very careful reading of the article}.
\\

\end{document}